\documentclass[aos,MSNbibl,nameyear,seceqn,dvips]{arximspdf}
\usepackage{mathrsfs,dcolumn}
\usepackage{graphicx}

\doi{10.1214/13-AOS1191} 
\volume{42}
\issue{1}
\pubyear{2014}
\firstpage{324}
\lastpage{351}

\makeatletter
\newcolumntype{d}[1]{D{.}{.}{#1}}
\newcommand{\rrvert}{\vert}
\newcommand{\llvert}{\vert}
\newcommand{\overset}{\stackrel}
\newtheorem{assumption}{Condition}
\newtheorem{lemma}{Lemma}
\newtheorem{proposition}{Proposition}
\newtheorem{theorem}{Theorem}
\newtheorem{corollary}{Corollary}
\makeatother

\begin{document}
\begin{frontmatter}

\title{Adaptive robust variable selection}
\runtitle{Adaptive robust variable selection}

\begin{aug}
\author[A]{\fnms{Jianqing} \snm{Fan}\thanksref{t2}\ead[label=e1]{jqfan@princeton.edu}},
\author[B]{\fnms{Yingying} \snm{Fan}\corref{}\thanksref{t3}\ead[label=e2]{fanyingy@marshall.usc.edu}}
\and
\author[C]{\fnms{Emre} \snm{Barut}\ead[label=e3]{abarut@us.ibm.com}}
\runauthor{J. Fan, Y. Fan and E. Barut}
\affiliation{Princeton University, University of Southern California
and\\ IBM T.~J. Watson Research Center}
\address[A]{J. Fan\\
Department of Operations Research\\
\quad and Financial Engineering\\
Princeton University\\
Princeton, New Jersey 08544\\
USA\\
\printead{e1}}
\address[B]{Y. Fan\\
Data Sciences and Operations Department\\
Marshall School of Business\\
University of Southern California\\
Los Angeles, California 90089\\
USA\\
\printead{e2}}
\address[C]{E. Barut\\
IBM T.~J. Watson Research Center\\
Yorktown Heights, New York 10598\\
USA \\
\printead{e3}}
\end{aug}
\thankstext{t2}{Supported by the National Institutes of Health Grant
R01-GM072611 and National Science Foundation Grants DMS-12-06464 and
DMS-07-04337.}
\thankstext{t3}{Supported by NSF CAREER Award DMS-1150318 and Grant
DMS-09-06784, the 2010 Zumberge Individual Award from USC's James H.
Zumberge Faculty
Research and Innovation Fund, and USC Marshall summer research funding.}

\received{\smonth{3} \syear{2013}}
\revised{\smonth{10} \syear{2013}}

%
\begin{abstract}
Heavy-tailed high-dimensional data are commonly encountered in various
scientific fields and pose great challenges to modern statistical
analysis. A~natural procedure to address this problem is to use
penalized \mbox{quantile} regression with weighted $L_1$-penalty, called
weighted robust Lasso \mbox{(WR-Lasso)}, in which weights are introduced to
ameliorate the bias problem induced by the $L_1$-penalty. In the
ultra-high dimensional setting, where the \mbox{dimensionality} can grow
exponentially with the sample size, we investigate the model selection
oracle property and establish the asymptotic normality of the WR-Lasso.
We show that only mild conditions on the model error distribution are
needed. Our theoretical results also reveal that adaptive choice of the
weight vector is essential for the WR-Lasso to enjoy these nice
asymptotic properties. To make the WR-Lasso practically feasible, we
propose a two-step procedure, called adaptive robust Lasso (AR-Lasso),
in which the weight vector in the second step is constructed based on
the $L_1$-penalized quantile regression estimate from the first step.
This two-step procedure is justified theoretically to possess the
oracle property and the asymptotic normality. Numerical studies
demonstrate the favorable finite-sample performance of the AR-Lasso.
\end{abstract}

%

\end{frontmatter}

\section{Introduction}
The advent of modern technology makes it easier to collect massive,
large-scale data-sets. A common feature of these data-sets is that the
number of covariates greatly exceeds the number of observations, a
regime opposite to conventional statistical settings. For example,
portfolio allocation with hundreds of stocks in finance involves a
covariance matrix of about tens of thousands of parameters, but the
sample sizes are often only in the order of hundreds (e.g., daily data
over a year period [\citet{fanfanlv08}]). Genome-wide association
studies in biology involve hundreds of thousands of single-nucleotide
polymorphisms (SNPs), but the available sample size is usually in
hundreds, also. Data-sets with large number of variables but relatively
small sample size pose great unprecedented challenges and opportunities
for statistical analysis.

Regularization methods have been widely used for high-dimensional
variable selection [\citet{tibshirani2}, \citet{fan2}, \citet{fan3},
\citet{Bickel06}, \citet{CandesTao07},
 \citet{BRT09}, \citet{lv1}, \citet{MB10}, \citet{Zhang07}, \citet{zou1}]. Yet, most
existing methods such as penalized least-squares or penalized
likelihood [\citet{lv2}] are designed for light-tailed
distributions. \citet{zhao2} established the irrepresentability
conditions for the model selection consistency of the Lasso estimator.
\citet{fan2} studied the oracle properties of nonconcave penalized
likelihood estimators for fixed dimensionality. \citet{lv1}
investigated the penalized least-squares estimator with folded-concave
penalty functions in the ultra-high dimensional setting and established
a nonasymptotic weak oracle property. \citet{FanLv2008} proposed and
investigated the sure independence screening method in the setting of
light-tailed distributions. The robustness of the aforementioned
methods have not yet been thoroughly studied and well understood.

Robust regularization methods such as the least absolute deviation
(LAD) regression and quantile regression have been used for variable
selection in the case of fixed dimensionality. See, for example, \citet{WLJ07},
\citet{LiY1}, \citet{zou3}, \citet{wuliu09}. The penalized composite likelihood
method was proposed in \citet{fanBradic11} for robust estimation in
ultra-high dimensions with focus on the efficiency of the method. They
still assumed sub-Gaussian tails. \citet{BC11} studied the
$L_1$-penalized quantile regression in high-dimensional sparse models
where the dimensionality could be larger than the sample size. We refer
to their method as robust Lasso (R-Lasso). They showed that the R-Lasso
estimate is consistent at the near-oracle rate, and gave conditions
under which the selected model includes the true model, and derived
bounds on the size of the selected model, uniformly in a compact set of
quantile indices. \citet{W2012} studied the $L_1$-penalized LAD
regression and showed that the estimate achieves near oracle risk
performance with a nearly universal penalty parameter and established
also a sure screening property for such an estimator. \citet{GM12}
obtained bounds on the prediction error of a large class of $L_1$-penalized estimators, including quantile regression. \citet{WWL2012}
considered the nonconvex penalized quantile regression in the
ultra-high dimensional setting and showed that the oracle estimate
belongs to the set of local minima of the nonconvex penalized quantile
regression, under mild assumptions on the error distribution. 

In this paper, we introduce the penalized quantile regression with the
weighted $L_1$-penalty
(WR-Lasso) for robust regularization, as in \citet{fanBradic11}. The
weights are introduced to reduce the bias problem induced by the
\mbox{$L_1$-}penalty. The flexibility of the choice of the weights provides
flexibility in shrinkage estimation of the regression coefficient.
WR-Lasso shares a similar spirit to the folded-concave penalized quantile-regression [\citet{ZL08}, \citet{WWL2012}], but avoids the
nonconvex optimization problem.
We establish conditions on the error distribution in order for the
WR-Lasso to successfully recover the true underlying sparse model with
asymptotic probability one. It turns out that the required condition is
much weaker than the sub-Gaussian assumption in \citet{fanBradic11}.
The only conditions we impose is that the density function of error has
Lipschitz property in a neighborhood around 0. This includes a large
class of heavy-tailed distributions such as the stable distributions,
including the Cauchy distribution. It also covers the double
exponential distribution whose density function is nondifferentiable at
the origin.

Unfortunately, because of the penalized nature of the estimator,
WR-Lasso estimate has a bias. In order to reduce the bias, the weights
in WR-Lasso need to be chosen adaptively according to the magnitudes of
the unknown true regression coefficients, which makes the bias
reduction infeasible for practical applications.

To make the bias reduction feasible, we introduce the adaptive robust
Lasso (AR-Lasso).
The AR-Lasso first runs R-Lasso to obtain an initial estimate, and then
computes the weight vector of the weighted $L_1$-penalty according to a
decreasing function of the magnitude of the initial estimate. After
that, AR-Lasso runs WR-Lasso with the computed weights. We formally
establish the model selection oracle property of AR-Lasso in the
context of \citet{fan2} with no assumptions made on the tail
distribution of the model error. In particular, the asymptotic
normality of the AR-Lasso is formally established.

This paper is organized as follows. First, we introduce our robust
estimators in Section~\ref{secnotation}. Then, to demonstrate the
advantages of our estimator, we show in Section~\ref{seclasso} with a
simple example that Lasso behaves suboptimally when noise has heavy
tails. In Section~\ref{secoracle}, we study the performance of the
oracle-assisted regularization estimator. Then in Section~\ref
{secselection}, we show that when the weights are adaptively chosen,
WR-Lasso has the model selection oracle property, and performs as well
as the oracle-assisted regularization estimate. In Section~\ref{secAN},
we prove the asymptotic normality of our proposed estimator.
The feasible estimator, AR-Lasso, is investigated in Section~\ref{secpenalty}.
Section~\ref{secsimul} presents the results of the simulation studies.
Finally, in Section~\ref{secproof}, we present the proofs of the main
theorems. Additional proofs, as well as the results of a genome-wide
association study, are provided in the supplementary Appendix
[\citet{ARLassoSupt}].

\section{Adaptive robust Lasso}\label{secnotation}
Consider the linear regression model
%
%
\begin{eqnarray}
\label{0} \mathbf{y}&=&\mathbf{X}\bolds{\beta}+ \bolds{\varepsilon},
\end{eqnarray}
where $\mathbf{y}$ is an $n$-dimensional response vector, $\mathbf
{X}=
(\mathbf{x}_1,\ldots,\mathbf{x}_n )^T= (\tilde\mathbf{x}_1,
\ldots, \tilde\mathbf{x}_p)$ is an $n\times p$ fixed design matrix,
$\bolds{\beta}= (\beta_1,\ldots,\beta_p )^T$ is a
$p$-dimensional regression coefficient vector, and $\bolds{\varepsilon
}=
(\varepsilon_1,\ldots,\varepsilon_n )^T$ is an $n$-dimensional
error vector whose components are independently distributed and satisfy
$P(\varepsilon_i\leq0) = \tau$ for some known constant $\tau\in
(0,1)$. Under this model, $\mathbf{x}_i^T \bolds{\beta}$ is the conditional
\mbox{$\tau$th-}quantile of $y_i$ given $\mathbf{x}_i$. We impose no
conditions on the heaviness of the tail probability or the
homoscedasticity of $\varepsilon_i$. We consider a challenging setting
in which $\log p = o(n^{b})$ with some constant $b > 0$. To ensure the
model identifiability and to enhance the model fitting accuracy and
interpretability, the true regression coefficient vector $\bolds{\beta
}^*$ is
commonly imposed to be sparse with only a small proportion of nonzeros
[\citet{tibshirani2}, \citet{fan2}]. Denoting the number of nonzero
elements of the true regression coefficients by $s_n$, we allow $s_n$
to slowly diverge with the sample size $n$ and assume that $s_n=o(n)$.
To ease the presentation, we suppress the dependence of $s_n$ on $n$
whenever there is no confusion. Without loss of generality, we write
$\bolds{\beta}^* = (\bolds{\beta}_{1}^{*T}, \mathbf{0}^T)^T$, that
is, only the
first $s$ entries are nonvanishing. The true model is denoted by
\[
\mathcal{M}_* = \operatorname{supp}\bigl(\bolds{\beta}^*\bigr)= \{
1,\ldots, s\}
\]
and its complement, $\mathcal{M}_*^c = \{s+1,\ldots, p\}$, represents the
set of noise variables.

We consider a fixed design matrix in this paper and denote by $\mathbf{S}
= (\mathbf{S}_1,\ldots,\mathbf{S}_n)^T = (\tilde\mathbf{x}_1,
\ldots, \tilde\mathbf{x}_s)$ the submatrix of $\mathbf{X}$
corresponding to the covariates
whose coefficients are nonvanishing. These variables will be referred
to as the signal covariates and the rest will be called noise
covariates. The set of columns that correspond to the noise covariates
is denoted by $\mathbf{Q}=(\mathbf{Q}_1,\ldots,\mathbf{Q}_{n})^T =
(\tilde\mathbf{x}_{s+1}, \ldots, \tilde\mathbf{x}_{p})$. We standardize each
column of $\mathbf{X}$ to have $L_2$-norm $\sqrt{n}$.

To recover the true model and estimate $\bolds{\beta}^{*}$, we
consider the
following regularization problem:
%
%
\begin{equation}
\label{genpenreg} \min_{\bolds{\beta}\in\mathbf{R}^p} \Biggl\{\sum
_{i=1}^n\rho_\tau\bigl(y_i
- \mathbf{x}_i^T\bolds{\beta}\bigr)+n
\lambda_n\sum_{j=1}^pp_{\lambda
_n}\bigl(|
\beta_j|\bigr) \Biggr\},
\end{equation}
where $\rho_\tau(u) = u(\tau- 1\{u\leq0\})$ is the quantile loss
function, and $p_{\lambda_n}(\cdot)$ is a nonnegative
penalty function on $[0, \infty)$ with a regularization parameter
$\lambda_n \geq0$. The use of quantile loss function in (\ref
{genpenreg}) is to overcome the difficulty of heavy tails of the error
distribution. Since $P(\varepsilon\leq0) = \tau$, (\ref{genpenreg})
can be interpreted as the sparse estimation of the conditional $\tau
$th quantile. Regarding the choice of $p_{\lambda_n}(\cdot)$, it was
demonstrated in \citet{lv1} and \citet{lv2} that folded-concave
penalties are more advantageous for variable
selection in high dimensions than the convex ones such as the
$L_1$-penalty. It is, however, computationally more challenging to
minimize the objective function in (\ref{genpenreg}) when $p_{\lambda
}(\cdot)$ is folded-concave. Noting that with a good initial estimate
$\hat{\bolds{\beta}}{}^{\mathrm{ini}}=(\hat\beta{}_1^{\mathrm{ini}}, \ldots,
\hat\beta{}_p^{\mathrm{ini}})^T$ of the true coefficient vector, we have
\[
p_{\lambda_n}\bigl(|\beta_j|\bigr)\approx p_{\lambda_n}\bigl(\bigl|\hat\beta
_j^{\mathrm{ini}}\bigr|\bigr) + p_{\lambda_n}'\bigl(\bigl|\hat
\beta_j^{\mathrm{ini}}\bigr|\bigr) \bigl( |\beta_j| - \bigl|\hat
\beta^{\mathrm{ini}}_j\bigr| \bigr).
\]
Thus, instead of (\ref{genpenreg}) we consider the following weighted
$L_1$-regularized quantile regression:
%
%
\begin{equation}
\label{eql1obj} L_n(\bolds{\beta})=\sum_{i=1}^n
\rho_\tau\bigl(y_i - \mathbf{x}_i^T
\bolds{\beta} \bigr)+n\lambda_n\|\mathbf{d}\circ\bolds{\beta}
\|_1,
\end{equation}
where $\mathbf{d}= (d_1,\ldots, d_p)^T$ is the vector of nonnegative
weights, and $\circ$ is the Hadamard product, that is, the
componentwise product of two vectors. This motivates us to define the
weighted robust Lasso (WR-Lasso) estimate as the global minimizer of
the convex function $L_n(\bolds{\beta})$ for a given nonstochastic
weight vector:
%
%
\begin{equation}
\label{eql1reg} \hat{\bolds{\beta}}= \mathop{\arg\min}_{\bolds{\beta
}}L_n(
\bolds{\beta}).
\end{equation}
The uniqueness of the global
minimizer is easily guaranteed by adding a negligible
$L_2$-regularization in implementation. In particular, when $d_j = 1$
for all $j$, the method will be referred to as robust Lasso (R-Lasso).

The\vspace*{1pt} adaptive robust Lasso (AR-Lasso) refers specifically to the
two-stage procedure in which the stochastic weights $\hat{d}_j =
p_{\lambda_n}'(|\hat\beta_j^{\mathrm{ini}}|)$ for $j=1, \ldots, p$ are
used in the second step for WR-Lasso and are constructed using a
concave penalty $p_{\lambda_n}(\cdot)$ and the initial estimates,
$\hat\beta_j^{\mathrm{ini}}$, from the first step. In practice, we
recommend using R-Lasso as the initial estimate and then using SCAD to
compute the weights in AR-Lasso.
The asymptotic result of this specific AR-Lasso is summarized in
Corollary~\ref{cor1} in Section~\ref{secpenalty} for the ultra-high dimensional robust
regression problem. This is a main contribution of the paper.

\section{Suboptimality of Lasso}\label{seclasso}
In this section, we use a specific example to illustrate that, in the
case of heavy-tailed error distribution, Lasso fails at model selection
unless the nonzero coefficients, $\beta^*_1,\ldots,\beta^*_s$, have
a very large magnitude. We assume that the errors $\varepsilon_1,
\ldots, \varepsilon_n$ have the identical symmetric stable
distribution and the characteristic function of $\varepsilon_1$ is
given by
\[
E\bigl[\exp(iu\varepsilon_1)\bigr] = \exp\bigl(-|u|^{\alpha}
\bigr),
\]
where $\alpha\in(0,2)$. By \citet{nolan12}, $E|\varepsilon_1|^p$ is
finite for $0<p<\alpha$, and $E|\varepsilon_1|^p = \infty$ for
$p\geq\alpha$. Furthermore, as $z \rightarrow\infty$,
\[
P\bigl(|\varepsilon_1| \geq z\bigr) \sim c_{\alpha}z^{-\alpha},
\]
where $c_{\alpha} = \sin(\frac{\pi\alpha}{2})\Gamma(\alpha)/\pi
$ is a constant depending only on $\alpha$, and we use the notation
$\sim$ to denote that two terms are equivalent up to some constant.
Moreover, for any constant vector $\mathbf{a}=(a_1, \ldots,
a_n)^T$, the linear combination $\mathbf{a}^T\bolds{\varepsilon}$
has the
following tail behavior:
%
%
\begin{equation}
\label{e040} P\bigl(\bigl|\mathbf{a}^T\bolds{\varepsilon}\bigr| > z\bigr) \sim
\|\mathbf{a}\| _{\alpha
}^{\alpha}c_{\alpha}z^{-\alpha}
\end{equation}
with $\|\cdot\|_{\alpha}$ denoting the $L_{\alpha}$-norm of a vector.

To demonstrate the suboptimality of Lasso, we consider a simple case in
which the design matrix satisfies the conditions that $\mathbf
{S}^T\mathbf{Q}=
\mathbf{0}$, $\frac{1}{n}\mathbf{S}^T\mathbf{S}= \mathbf{I}_s$,
the columns of $\mathbf{Q}$
satisfy $|\operatorname{supp}(\tilde\mathbf{x}_j)|=m_n = O(n^{1/2})$
and $\operatorname{supp}
(\tilde\mathbf{x}_k) \cap\operatorname{supp}(\tilde\mathbf{x}_j)
= \varnothing
$ for any $k\neq j$ and $k, j\in\{s+1, \ldots, p\}$. Here, $m_n$ is a
positive integer measuring the sparsity level of the columns of
$\mathbf{Q}$.
We assume that there are only fixed number of true variables, that is,
$s$ is finite, and that $\max_{ij} | x_{ij} |=O(n^{1/4})$.
Thus, it is easy to see that $p = O(n^{1/2})$. In addition, we assume
further that all nonzero regression coefficients are the same and
$\beta^*_1 = \cdots= \beta^*_s = \beta_0>0$.

We first consider R-Lasso, which is the global minimizer of (\ref
{eql1reg}). We will later see in Theorem~\ref{T2} that by choosing
the tuning parameter
\[
\lambda_n=O\bigl((\log n)^{2}\sqrt{(\log p)/n}\bigr),
\]
R-Lasso can recover the true support $\mathcal{M}_*=\{1, \ldots, s\}$
with probability tending to~1. Moreover, the signs of the true
regression coefficients can also be recovered with asymptotic
probability one as long as the following condition on signal strength
is satisfied:
%
%
\begin{equation}
\label{eql1sigstrength} \lambda_n^{-1} \beta_0 \to
\infty,\qquad\mbox{that is, }(\log n)^{-2}\sqrt{n/(\log p)}
\beta_0 \rightarrow\infty.
\end{equation}

Now, consider Lasso, which minimizes
%
%
\begin{equation}
\label{eqlasso} \widetilde L_n(\bolds{\beta}) = \tfrac{1}{2}\|
\mathbf{y}- \mathbf{X}\bolds{\beta}\|_2^2 +n\lambda
_n\|\bolds{\beta}\|_1.
\end{equation}
We will see that for (\ref{eqlasso}) to recover the true model and
the correct signs of coefficients, we need a much stronger signal level
than that is given in (\ref{eql1sigstrength}). By results in
optimization theory, the Karush--Kuhn--Tucker (KKT) conditions
guaranteeing the necessary and sufficient conditions for $\tilde{\bolds
{\beta}}$ with $\mathcal{M} = \operatorname{supp}(\tilde{\bolds{\beta
}})$ being a
minimizer to~(\ref{eqlasso}) are
\begin{eqnarray*}
\tilde{\bolds{\beta}}_{\mathcal{M}}+ n\lambda_n\bigl(\mathbf
{X}_{\mathcal{M}}^T\mathbf{X} _{\mathcal{M}}\bigr)^{-1}
\operatorname{sgn}(\tilde{\bolds{\beta}}_{\mathcal{M}}) &=& \bigl(\mathbf{X}
_{\mathcal{M}}^T\mathbf{X}_{\mathcal{M}}\bigr)^{-1}
\mathbf{X}_{\mathcal
{M}}^T\mathbf{y},
\\
\bigl\|\mathbf{X}_{\mathcal{M}^c}^T(\mathbf{y}-\mathbf{X}_{\mathcal
{M}}
\tilde{\bolds{\beta}} _{\mathcal{M}})\bigr\|_\infty&\leq& n
\lambda_n,
\end{eqnarray*}
where $\mathcal{M}^c$ is the complement of $\mathcal{M}$, $\bolds
{\beta}
_{\mathcal{M}}$ is the subvector formed by entries of~$\bolds{\beta
}$ with
indices in $\mathcal{M}$, and $\mathbf{X}_{\mathcal{M}}$ and
$\mathbf{X}
_{\mathcal{M}^c}$ are the submatrices formed by columns of $\mathbf
{X}$ with
indices in $\mathcal{M}$ and $\mathcal{M}^c$, respectively. It is
easy to see from the above two conditions that for Lasso to enjoy the
sign consistency, $\operatorname{sgn}(\tilde{\bolds{\beta}}) =
\operatorname{sgn}(\bolds{\beta}^*)$ with
asymptotic probability one, 
we must have these two conditions satisfied with $\mathcal{M}=\mathcal
{M}^*$ with probability tending to 1. Since we have assumed that
$\mathbf{Q}
^T\mathbf{S}= \mathbf{0}$ and $n^{-1}\mathbf{S}^T\mathbf{S}=\mathbf
{I}$, the above sufficient and
necessary conditions can also be written as
%
%
\begin{eqnarray}
\tilde{\bolds{\beta}}_{\mathcal{M}^*}+ \lambda_n\operatorname{sgn}(
\tilde{\bolds{\beta}} _{\mathcal{M}^*}) &=& \bolds{\beta}_{\mathcal{M}^*}^*+
n^{-1}\mathbf{S}^T\bolds{\varepsilon},\label{eqkkt1a}
\\
\bigl\|\mathbf{Q}^T\bolds{\varepsilon}\bigr\|_\infty&\leq& n
\lambda_n. \label{eqkkt2a}
\end{eqnarray}
Conditions (\ref{eqkkt1a}) and (\ref{eqkkt2a}) are hard for Lasso
to hold simultaneously. The following proposition summarizes the
necessary condition, whose proof is given in the supplementary material
[\citet{ARLassoSupt}].

%
\begin{proposition}\label{prop1}
In the above model, with probability at least $1-e^{-\tilde c_0}$,
where $\tilde c_0$ is some positive constant, Lasso does not have sign
consistency, unless the following signal condition is satisfied
%
%
\begin{equation}
\label{eqlassosigstrength} n^{(3/4)-(1/\alpha)}\beta
_0\rightarrow\infty.
\end{equation}
\end{proposition}

Comparing this with (\ref{eql1sigstrength}), it is easy to see that
even in this simple case, Lasso needs much stronger signal levels than
R-Lasso in order to have a sign consistency in the presence of a
heavy-tailed distribution.

\section{Model selection oracle property}\label{secmain}

In this section, we establish the model selection oracle property of
WR-Lasso. The study enables us to see the bias due to penalization, and
that an adaptive weighting scheme is needed in order to eliminate such
a bias. We need the following condition on the distribution of noise.

%
\begin{assumption}\label{assp1}
There exist universal constants $c_1>0$ and $c_2>0$ such that for any
$u$ satisfying $ |u|\leq c_1$, $f_i(u)$'s are uniformly bounded away
from 0 and $\infty$ and
\[
\bigl|F_i(u)-F_i(0) - uf_i(0) \bigr|\leq
c_2u^2,
\]
where $f_i(u)$ and $F_i(u)$ are the density function and distribution
function of the error $\varepsilon_i$, respectively.
\end{assumption}

Condition~\ref{assp1} implies basically that each $f_i(u)$ is Lipschitz around
the origin. Commonly used distributions such as the double-exponential
distribution and stable distributions including the Cauchy distribution
all satisfy this condition.

Denote by $\mathbf{H}= \operatorname{diag}\{f_1(0),\ldots, f_n(0)\}$.
The next
condition is on the submatrix of $\mathbf{X}$ that corresponds to signal
covariates and the magnitude of the entries of $\mathbf{X}$.

%
\begin{assumption}\label{assp2}
The\vspace*{1pt} eigenvalues of $\frac{1}{n}\mathbf{S}^T\mathbf{H}\mathbf{S}$
are bounded from below
and above by some positive constants $c_0$ and $1/c_0$, respectively.
Furthermore,
\[
\kappa_n\equiv\max_{ij}|x_{ij}| = o
\bigl(\sqrt{n}s^{-1}\bigr).
\]
\end{assumption}
Although Condition~\ref{assp2} is on the fixed design matrix, we note
that the above condition on $\kappa_n$ is satisfied with asymptotic
probability one when the design matrix is generated from some
distributions. For instance, if the entries of $\mathbf{X}$ are independent
copies from a subexponential distribution, the bound on $\kappa_n$ is
satisfied with asymptotic probability one as long as $s=o (\sqrt
{n}/(\log p) )$; if the components are generated from sub-Gaussian
distribution, then the condition on $\kappa_n$ is satisfied with
probability tending to one when $s=o (\sqrt{n/(\log p)} )$.

\subsection{Oracle regularized estimator}\label{secoracle}
To evaluate our newly proposed method, we first study how well one can
do with the assistance of the oracle information on the locations of
signal covariates. Then we use this to establish the asymptotic
property of our estimator without the oracle assistance. Denote by
$\hat{\bolds{\beta}}{}^o = ((\hat{\bolds{\beta}}{}_1^o)^T,
\mathbf{0}^T)^T$ the
oracle regularized estimator (ORE) with $\hat{\bolds{\beta}}{}_1^o
\in
\mathbf{R}^s$ and $\mathbf{0}$ being the vector of all zeros, which
minimizes $L_n(\bolds{\beta})$ over the space $\{\bolds{\beta
}=(\bolds{\beta}_1^T, \bolds{\beta}
_2^T)^T\in\mathbf{R}^p\dvtx  \bolds{\beta}_2 = \mathbf{0}\in\mathbf
{R}^{p-s} \}$.
The next theorem shows that ORE is consistent, and estimates the
correct sign of the true coefficient vector with probability tending to
one. We use $\mathbf{d}_0$ to denote the first $s$ elements of~$\mathbf{d}$.

%
\begin{theorem} \label{T1}
Let $\gamma_n = C_1(\sqrt{s(\log n)/n}+ \lambda_n\|\mathbf{d}_0\|
_2)$ with
$C_1>0$ a constant. If Conditions~\ref{assp1} and~\ref{assp2} hold
and $\lambda_n\|\mathbf{d}_0\|_2\sqrt{s}\kappa_n\rightarrow0$,
then there
exists some constant $c>0$ such that
%
%
\begin{equation}
P \bigl(\bigl\|\hat{\bolds{\beta}}{}_1^o - \bolds{
\beta}^*_1\bigr\|_2 \leq\gamma_n \bigr)\geq1-
n^{-cs}.
\end{equation}
If in addition $\gamma_n^{-1}\min_{1\leq j\leq s}|\beta_j^*|
\rightarrow\infty$, then with probability at least $1- n^{-c s}$,
\[
\operatorname{sgn}\bigl(\hat{\bolds{\beta}}{}_1^o\bigr) =
\operatorname{sgn}\bigl(\bolds{\beta}_1^*\bigr),
\]
where the above equation should be understood componentwisely.
\end{theorem}

As shown in Theorem~\ref{T1}, the consistency rate of $\hat{\bolds{\beta}}{}_1^o$ in terms of the vector \mbox{$L_2$-}norm is given by $\gamma_n$. The
first component of $\gamma_n$, $C_1\sqrt{s(\log n)/n}$, is the oracle
rate within a factor of $\log n$, and the second component $C_1\lambda
_n\|\mathbf{d}_0\|_2$ reflects the bias due to penalization. If no prior
information is available, one may choose equal weights $\mathbf{d}_0 = (1,
1,\ldots, 1)^T$, which corresponds to R-Lasso. Thus, for R-Lasso, with
probability at least $1-n^{-cs}$, it holds that
%
%
\begin{equation}
\label{eq12a} \bigl\|\hat{\bolds{\beta}}{}_1^o - \bolds{
\beta}_1^*\bigr\|_2 \leq C_1\bigl(\sqrt{s(\log
n)/n}+\sqrt{s}\lambda_n\bigr).
\end{equation}

\subsection{WR-Lasso}\label{secselection}

In this section, we show that even without the oracle information,
WR-Lasso enjoys the same asymptotic property as in Theorem~\ref{T1}
when the weight vector is appropriately chosen. Since the regularized
estimator $\hat{\bolds{\beta}}$ in (\ref{eql1reg}) depends on
the full
design matrix $\mathbf{X}$, we need to impose the following conditions
on the
design matrix to control the correlation of columns in $\mathbf{Q}$
and $\mathbf{S}$.

%
\begin{assumption}\label{assp3}
With $\gamma_n$ defined in Theorem~\ref{T1}, it holds that
\[
\biggl\|\frac{1}{n}\mathbf{Q}^T\mathbf{H}\mathbf{S}
\biggr\|_{2,\infty
} < \frac{\lambda
_n}{2\|\mathbf{d}_1^{-1}\|_\infty\gamma_n},
\]
where $\|\mathbf{A}\|_{2,\infty} = \sup_{\mathbf{x}\neq0}{\|
\mathbf{A}\mathbf
{x}\|_\infty}/{\|\mathbf{x}\|_2}$ for a matrix $\mathbf{A}$ and vector
$\mathbf{x}$, and $\mathbf{d}_1^{-1}=(d_{s+1}^{-1},\ldots, d_p^{-1})^T$.
Furthermore, $\log(p)= o(n^b)$ for some constant $b \in(0,1)$.
\end{assumption}

To understand the implications of Condition~\ref{assp3}, we consider
the case of 
$f_1(0)=\cdots= f_n(0)\equiv f(0)$. In the special case of $\mathbf
{Q}^T\mathbf{S}
= \mathbf{0}$, Condition~\ref{assp3} is satisfied automatically. In the
case of equal correlation, that is, $n^{-1}\mathbf{X}^T\mathbf{X}$ having
off-diagonal elements all equal to $\rho$, the above Condition~\ref{assp3} reduces to 
\[
|\rho| < \frac{\lambda_n}{4f(0)\|\mathbf{d}_1^{-1}\|_\infty\sqrt
{s}\gamma_n}.
\]
This puts an upper bound on the correlation coefficient $\rho$ for
such a dense matrix.

It is well known that for Gaussian errors, the optimal choice of
regularization parameter $\lambda_n$ has the order $\sqrt{(\log
p)/n}$ [\citet{BRT09}]. The distribution of the model noise with
heavy tails demands a larger choice of $\lambda_n$ to filter the noise
for R-Lasso. When $\lambda_n \geq\sqrt{(\log n) /n}$, $\gamma_n$
given in (\ref{eq12a}) is in the order of $C_1\lambda_n\sqrt{s}$. In
this case, Condition~\ref{assp3} reduces to
%
%
\begin{equation}
\label{eq13} \bigl\|n^{-1}\mathbf{Q}^T\mathbf{H}\mathbf{S}
\bigr\|_{2,\infty} < O\bigl(s^{-1/2}\bigr).
\end{equation}
For WR-Lasso, if the weights are chosen such that $\|\mathbf{d}_0\|
_2=O(\sqrt{s(\log n)/n}/\lambda_n)$ and $\|\mathbf{d}_1\|_{\infty
}=O(1)$, then
$\gamma_n$ is in the order of $C_1\sqrt{s(\log n)/n}$, and
correspondingly, Condition~\ref{assp3} becomes
\[
\bigl\|n^{-1}\mathbf{Q}^T\mathbf{H}\mathbf{S}\bigr\|_{2,\infty}
< O \bigl(\lambda_n\sqrt{n/\bigl(s(\log n)\bigr)} \bigr).
\]
This is a more relaxed condition than (\ref{eq13}), since with
heavy-tailed errors, the optimal $\lambda_n$ should be larger than
$\sqrt{(\log p)/n}$. In other words, WR-Lasso not only reduces the
bias of the estimate, but also allows for stronger correlations among
the signal and noise covariates. However, the above choice of weights
depends on unknown locations of signals. A data-driven choice will be given
in Section~\ref{secpenalty}, in which the resulting AR-Lasso
estimator will be studied.

The following theorem shows the model selection oracle property of the
WR-Lasso estimator.

%
\begin{theorem}\label{T2}
Suppose Conditions~\ref{assp1}--\ref{assp3} hold. In addition,
assume that $\min_{j\geq s+1}\,d_j >c_3$ with some constant $c_3>0$,
%
%
\begin{equation}
\label{eq12} \gamma_ns^{3/2}\kappa_n^2(
\log_2n)^2 = o\bigl(n \lambda_n^2
\bigr),\qquad\lambda_n\|\mathbf{d}_0\|_2
\kappa_n\max\bigl\{\sqrt{s},\| \mathbf{d}_0\|_2\bigr\}
\rightarrow0
\end{equation}
and $\lambda_n>2\sqrt{(1+c)(\log p)/n}$, where $\kappa_n$ is defined
in Condition~\ref{assp2}, $\gamma_n$ is defined in Theorem~\ref{T1},
and $c$ is some positive constant. Then, with probability at least
$1-O(n^{-cs})$, there exists a global minimizer $\hat{\bolds{\beta}}
=((\hat{\bolds{\beta}}{}_1^o)^T, \hat{\bolds{\beta}}{}_2^T)^T$ of
$L_n(\bolds{\beta})$
which satisfies
\begin{longlist}[(2)]
\item[(1)] $\hat{\bolds{\beta}}_2 = 0$;
\item[(2)] $\|\hat{\bolds{\beta}}{}_1^o - \bolds{\beta}_1^*\|_2
\leq\gamma_n$.
\end{longlist}
\end{theorem}

Theorem~\ref{T2} shows that the WR-Lasso estimator enjoys the same
property as ORE with probability tending to one. However, we impose
nonadaptive assumptions on the weight vector $\mathbf{d}= (\mathbf
{d}_0^T, \mathbf{d}
_1^T)^T$. For noise covariates, we assume $\min_{j>s}\,d_j>c_3$, which
implies that each coordinate needs to be penalized. For the signal
covariates, we impose (\ref{eq12}), which requires $\|\mathbf{d}_0\|
_2$ to
be small.

When studying the nonconvex penalized quantile regression, \citet{WWL2012} assumed that $\kappa_n$ is bounded and the density functions
of $\varepsilon_i$'s are uniformly bounded away from 0 and $\infty$
in a small neighborhood of 0. Their assumption on the error
distribution is weaker than our Condition~\ref{assp1}. We remark that
the difference is because we have weaker conditions on $\kappa_n$ and
the penalty function [see Condition~\ref{assp2} and (\ref{eq12})]. In
fact, our Condition~\ref{assp1} can be weakened to the same condition
as that in \citet{WWL2012} at the cost of imposing stronger assumptions
on $\kappa_n$ and the weight vector~$\mathbf{d}$.

\citet{BC11} and \citet{W2012} imposed the restricted eigenvalue
assumption of the design matrix and studied the $L_1$-penalized\vadjust{\goodbreak}
quantile regression and LAD regression, respectively. We impose
different conditions on the design matrix and allow flexible shrinkage
by choosing $\mathbf{d}$. In addition, our Theorem~\ref{T2} provides a
stronger result than consistency; we establish model selection oracle
property of the estimator.

\subsection{Asymptotic normality}\label{secAN}
We now present the asymptotic normality of our estimator. Define
$\mathbf{V}
_n = (\mathbf{S}^T\mathbf{H}\mathbf{S})^{-1/2}$ and $\mathbf{Z}_n =
(\mathbf{Z}_{n1}, \ldots, \mathbf{Z}
_{nn})^T = \mathbf{S}\mathbf{V}_n$ with \mbox{$\mathbf{Z}_{nj}\in\mathbf
{R}^s$} for
$j=1,\ldots, n$. 

%
\begin{theorem}\label{T3} Assume the conditions of Theorem~\ref{T2}
hold, the first and second order derivatives $f_i'(u)$ and $f_i''(u)$
are uniformly bounded in a small neighborhood\vspace*{1pt} of~0 for all $i=1,\ldots,
n$, and that $\|\mathbf{d}_0\|_2 = O(\sqrt{s/n} /\lambda_n)$,\break
$\max_i\|\mathbf{H}^{1/2}\mathbf{Z}_{ni}\|_2 = o(s^{-7/2}(\log
s)^{-1})$, and
$\sqrt{n/s}\min_{1\leq j\leq s}|\beta_j^*| \rightarrow\infty$.
Then, with probability tending to 1 there exists a global minimizer
$\hat{\bolds{\beta}}= ((\hat{\bolds{\beta}}{}_1^o)^T, \hat{\bolds{\beta}}{}_2^T)^T$
of $L_n(\bolds{\beta})$ such that $\hat{\bolds{\beta}}_2 = 0$. Moreover,
\[
\mathbf{c}^T\bigl(\mathbf{Z}_n^T
\mathbf{Z}_n\bigr)^{-1/2}\mathbf{V}_n^{-1}
\biggl[\bigl(\hat{\bolds{\beta}}{}_1^o - \bolds{
\beta}_1^*\bigr) + \frac{n\lambda_n}{2}\mathbf{V}_n^{2}
\tilde\mathbf{d}_0 \biggr] \overset{\mathscr{D}} {\longrightarrow}N
\bigl(0, \tau(1-\tau) \bigr),
\]
where $\mathbf{c}$ is an arbitrary $s$-dimensional vector satisfying
$\mathbf{c}
^T\mathbf{c}=1$, and $\tilde\mathbf{d}_0$ is an \mbox{$s$-}dimensional
vector with the
$j$th element $d_j\operatorname{sgn}(\beta^*_j)$.
\end{theorem}

The proof of Theorem~\ref{T3} is an extension of the proof on the
asymptotic normality theorem for the LAD estimator in \citet{Pollard1990}, in which the theorem is proved for fixed dimensionality.
The idea is to approximate $L_n(\bolds{\beta}_1,\mathbf{0})$ in
(\ref
{eql1reg}) by a sequence of quadratic functions, whose minimizers
converge to normal distribution. Since $L_n(\bolds{\beta}_1, \mathbf
{0})$ and
the quadratic approximation are close, their minimizers are also close,
which results in the asymptotic normality in Theorem~\ref{T3}.\vspace*{1pt}

Theorem~\ref{T3} assumes that $\max_i\|\mathbf{H}^{1/2}\mathbf
{Z}_{ni}\|_2 =
o(s^{-7/2}(\log s)^{-1})$. Since by definition $\sum_{i=1}^n\|\mathbf{H}
^{1/2}\mathbf{Z}_{ni}\|_2^2 = s$, it is seen that the condition
implies $s =
o(n^{1/8})$. This assumption is made to guarantee that the quadratic
approximation is close enough to $L_n(\bolds{\beta}_1, \mathbf{0})$.
When $s$ is
finite, the condition becomes $\max_i\|\mathbf{Z}_{ni}\|_2 = o(1)$,
as in
\citet{Pollard1990}. Another important assumption is $\lambda_n\sqrt
{n}\|\mathbf{d}_0\|_2 = O(\sqrt{s})$, which is imposed to make sure
that the
bias $2^{-1}n\lambda_n\mathbf{c}^T\mathbf{V}_n \tilde\mathbf{d}_0$
caused by the
penalty term does not diverge. For instance, using R-Lasso will create
a nondiminishing bias, and thus cannot be guaranteed to have asymptotic
normality.

Note that we do not assume a parametric form of the error distribution.
Thus, our oracle estimator is in fact a semiparametric estimator with
the error density as the nuisance parameter. Heuristically speaking,
Theorem~\ref{T3} shows that the asymptotic variance of $\sqrt{n}(\hat
{\bolds{\beta}}{}_1^o - \bolds{\beta}_1^*)$ is $n\tau
(1-\tau)\mathbf{V}_n\mathbf{Z}
_n^T\mathbf{Z}_n\mathbf{V}_n$. Since $\mathbf{V}_n = (\mathbf
{S}^T\mathbf{H}\mathbf{S})^{-1/2}$ and $\mathbf{Z}_n
= \mathbf{S}\mathbf{V}_n$, if the model errors $\varepsilon_i$ are
i.i.d. with
density function $f_{\varepsilon}(\cdot)$, then this asymptotic
variance reduces to $\tau(1-\tau)(n^{-1}f_{\varepsilon}^2(0)\mathbf{S}
^T\mathbf{S})^{-1}$. In the random design case where the true covariate
vectors $\{\mathbf{S}_i\}_{i=1}^n$ are i.i.d. observations,
$n^{-1}\mathbf{S}^T\mathbf{S}
$ converges to $E[\mathbf{S}_1^T\mathbf{S}_1]$ as $n \rightarrow
\infty$, and the
asymptotic variance reduces to $\tau(1-\tau)(f_{\varepsilon
}^2(0)E[\mathbf{S}_1^T\mathbf{S}_1])^{-1}$. This is the
semiparametric efficiency
bound derived by \citet{NP1990} for random designs.
In fact, if we assume that $(\mathbf{x}_i, y_i)$ are i.i.d., then the
conditions of Theorem~\ref{T3} can hold with asymptotic probability one.
Using similar arguments, it can be formally shown that $\sqrt{n}(\hat
{\bolds{\beta}}{}_1^o - \bolds{\beta}_1^*)$ is
asymptotically normal with
covariance matrix equal to the aforementioned semiparametric efficiency
bound. Hence, our oracle estimator is semiparametric efficient.


\section{Properties of the adaptive robust Lasso} \label{secpenalty}

In previous sections, we have seen that the choice of the weight vector
$\mathbf{d}$ plays a pivotal role for the WR-Lasso estimate to enjoy the
model selection oracle property and asymptotic normality. In fact,
conditions in Theorem~\ref{T2} require that $\min_{j \ge s+1}\,d_j
>c_3$ and that $\|\mathbf{d}_0\|_2$ does not diverge too fast. Theorem
\ref
{T3} imposes an even more stringent condition, $\|\mathbf{d}_0\|
_2=O(\sqrt{s/n} / \lambda_n )$, on the weight vector $\mathbf{d}_0$.
For R-Lasso,
$\|\mathbf{d}_0\|_2=\sqrt{s}$ and these conditions become very restrictive.
For example, the condition in Theorem~\ref{T3} becomes $\lambda_n =
O(n^{-1/2})$, which is too low for a thresholding level even for
Gaussian errors. Hence, an adaptive choice of weights is needed to
ensure that those conditions are satisfied. To this end, we propose a
two-step procedure.

In the first step, we use R-Lasso, which gives the estimate $\hat{\bolds
{\beta}}{}^{\mathrm{ini}}$. As has been shown in \citet{BC11} and \citet{W2012},
R-Lasso is consistent at a near-oracle rate $\sqrt{s(\log p)/n}$ and
selects the true model $\mathcal{M}^*$ as a submodel [in other words,
R-Lasso has the sure screening property using the terminology of \citet{FanLv2008}] with asymptotic probability one, namely,
\[
\operatorname{supp} \bigl(\hat{\bolds{\beta}}{}^{\mathrm{ini}} \bigr) \supseteq
\operatorname{supp} \bigl(\bolds{\beta} ^* \bigr)\quad\mbox{and}\quad\bigl\|
\hat{
\bolds{\beta}}{}^{\mathrm{ini}}_1 - \bolds{\beta}^*_1
\bigr\|_2 = O\bigl(\sqrt{s(\log p)/n}\bigr).
\]
We remark that our Theorem~\ref{T2} also ensures the consistency of
R-Lasso. Compared to \citet{BC11}, Theorem~\ref{T2} presents stronger
results but also needs more restrictive conditions for R-Lasso. As will
be shown in latter theorems, only the consistency of R-Lasso is needed
in the study of AR-Lasso, so we quote the results and conditions on
R-Lasso in \citet{BC11} with the mind of imposing weaker conditions.

In the second step, we set $\hat{\mathbf{d}}= (\hat d_1,\ldots,
\hat d_p)^T$ with $\hat d_j = p'_{\lambda_n}(|\hat\beta^{\mathrm{ini}}_j|)$
where $p_{\lambda_n}(|\cdot|)$ is a folded-concave penalty function,
and then solve the regularization problem (\ref{eql1reg}) with a
newly computed weight vector.
Thus, vector $\hat{\mathbf{d}}_0$ is expected to be close to the vector
$(p_{\lambda_n}'(|\beta^*_1|), \ldots, p'_{\lambda_n}(|\beta
_s^*|))^T$ under $L_2$-norm. If a folded-concave penalty such as SCAD
is used, then $p'_{\lambda_n}(|\beta_j^*|)$ will be close, or even
equal, to zero for $1\leq j\leq s$, and thus the magnitude of $\|
\hat{\mathbf{d}}_0\|_2$ is negligible.\vadjust{\goodbreak}

Now, we formally establish the asymptotic properties of AR-Lasso. We
first present a more general result and then highlight our recommended
procedure, which uses R-Lasso as the initial estimate and then uses
SCAD to compute the stochastic weights, in Corollary~\ref{cor1}. Denote\vspace*{-2pt} by %
$\mathbf{d}^*=(d_1^*,\ldots, d_p^*)$ with $d_j^* =p'_{\lambda
_n}(|\beta
_j^*|)$. Using the weight vector $\hat{\mathbf{d}}$, AR-Lasso minimizes
the following objective function:
%
%
\begin{equation}
\label{eql1reghat} \widehat L_n(\bolds{\beta}) = \sum
_{i=1}^n\rho_\tau\bigl(y_i -
\mathbf{x}_i^T\bolds{\beta}\bigr) +n
\lambda_n\|\hat{\mathbf{d}}\circ\bolds{\beta}\|_1.
\end{equation}
We also need the following conditions to show the model selection
oracle property of the two-step procedure.

%
\begin{assumption}\label{assp5}
With asymptotic probability one, the initial estimate satisfies $\|
\hat{\bolds{\beta}}{}^{\mathrm{ini}}-\bolds{\beta}^*\|_2\leq C_2\sqrt{s(\log
p)/n}$ with
some constant $C_2>0$.
\end{assumption}
As discussed above, if R-Lasso is used to obtain the initial estimate,
it satisfies the above condition. Our second condition is on the
penalty function.

%
\begin{assumption}\label{assp4} $p_{\lambda_n}'(t)$ is nonincreasing
in $t\in(0,\infty)$ and is Lipschitz with constant $c_5$, that is,
\[
\bigl|p'_{\lambda_n}\bigl(|\beta_1|\bigr) - p'_{\lambda_n}\bigl(|
\beta_2|\bigr)\bigr|\leq c_5|\beta_1 -
\beta_2|
\]
for any $\beta_1, \beta_2 \in\mathbf{R}$. Moreover, $p'_{\lambda
_n}(C_2\sqrt{s(\log p)/n}) > \frac{1}{2} p'_{\lambda_n}(0+)$ for
large enough~$n$, where $C_2$ is defined in Condition~\ref{assp5}.
\end{assumption}

For the SCAD [\citet{fan2}] penalty, $p_{\lambda_n}'(\beta)$ is
given by
%
%
\begin{equation}
\label{eqscad} p'_{\lambda_n}(\beta)= 1\{\beta\leq
\lambda_n\}+ \frac{ (a
\lambda_n - \beta)_{+}}{(a-1)\lambda_n}1\{\beta> \lambda_n\}
\end{equation}
for a given constant $a > 2$, and it can be easily verified that
Condition~\ref{assp4} holds if $\lambda_n>2(a+1)^{-1}C_2\sqrt{s(\log p)/n}$.

%
\begin{theorem}\label{T4} Assume conditions of Theorem~\ref{T2} hold
with $\mathbf{d}= \mathbf{d}^*$ and $\gamma_n = a_n$, where
\[
a_n = C_3 \bigl(\sqrt{s(\log n)/n}+
\lambda_n \bigl(\bigl\|\mathbf{d}_0^*\bigr\| _2 +
C_2 c_5\sqrt{s(\log p)/n} \bigr) \bigr)
\]
with some constant $C_3>0$ and $\lambda_ns\kappa_n\sqrt{(\log p)/n}
\rightarrow0$.
Then, under Conditions~\ref{assp5} and~\ref{assp4}, with probability
tending to one, there exists a global minimizer $\hat{\bolds{\beta}}=
(\hat{\bolds{\beta}}{}_1^T, \hat{\bolds{\beta}}{}_2^T)^T$ of
(\ref{eql1reghat})
such that $\hat{\bolds{\beta}}_2 = \mathbf{0}$ and
$\|\hat{\bolds{\beta}}_1 - \bolds{\beta}_1^*\|_2 \leq a_n$.
\end{theorem}

The results in Theorem~\ref{T4} are analogous to those in Theorem~\ref{T2}. The extra term $\lambda_n\sqrt{s(\log p)/n}$ in the convergence
rate $a_n$, compared to the convergence rate $\gamma_n$ in Theorem
\ref{T2}, is caused by the bias of the initial estimate $\hat{\bolds
{\beta}}{}^{\mathrm{ini}}$. Since the regularization parameter $\lambda_n$
goes to
zero, the bias of AR-Lasso is much smaller than that of the initial
estimator $\hat{\bolds{\beta}}{}^{\mathrm{ini}}$. Moreover, the AR-Lasso
$\hat{\bolds{\beta}}$ possesses the model selection oracle property.

Now we present the asymptotic normality of the AR-Lasso estimate.

%
\begin{assumption}\label{assp6} The smallest signal satisfies $\min
_{1\leq j\leq s}|\beta_j^*| >\break  2C_2\sqrt{(s\log p)/n}$. Moreover, it
holds that
$p''_{\lambda_n}(|\beta|)=o(s^{-1}\lambda_n^{-1}(n\log p)^{-1/2})$
for any $|\beta|> 2^{-1}\min_{1\leq j\leq s}|\beta_j^*|$.
\end{assumption}

The above condition on the penalty function is satisfied when the SCAD
penalty is used and $\min_{1\leq j\leq s}|\beta_j^*|\geq2a\lambda
_n$ where $a$ is the parameter in the SCAD penalty~(\ref{eqscad}).

%
\begin{theorem}\label{T5}
Assume conditions of Theorem~\ref{T3} hold with $\mathbf{d}=\mathbf
{d}^*$ and
$\gamma_n = a_n$, where $a_n$ is defined in Theorem~\ref{T4}. Then,
under Conditions~\ref{assp5}--\ref{assp6}, with asymptotic
probability one, there exists a global minimizer $\hat{\bolds{\beta}}$ of
(\ref{eql1reghat}) having the same asymptotic properties as those in
Theorem~\ref{T3}.
\end{theorem}

With the SCAD penalty, conditions in Theorems~\ref{T4} and~\ref{T5}
can be simplified and AR-Lasso still enjoys the same asymptotic
properties, as presented in the following corollary.

%
\begin{corollary}\label{cor1} Assume $\lambda_n = O(\sqrt{s(\log
p)(\log\log n)/n})$, $\log p = o(\sqrt{n})$, $\min_{1\leq j\leq
s}|\beta_j^*|\geq2a\lambda_n$ with $a$ the parameter in the SCAD
penalty and $\kappa_n=o(n^{1/4}s^{-1/2}(\log n)^{-3/2}(\log
p)^{1/2})$. Further assume that $\|n^{-1}\mathbf{Q}^T\mathbf
{H}\mathbf{S}\|_{2,\infty} <\break
C_4\sqrt{(\log p)(\log\log n)/\log n}$ with $C_4$ some positive
constant. Then, under Conditions~\ref{assp1} and~\ref{assp2}, with
asymptotic probability one, there exists a global minimizer $\hat{\bolds
{\beta}}= (\hat{\bolds{\beta}}{}_1^T, \hat{\bolds{\beta}}{}_2^T)^T$ of
$\widehat L_n(\bolds{\beta})$ such that
\[
\bigl\|\hat{\bolds{\beta}}_1 - \bolds{\beta}_1^*
\bigr\|_2 \leq O\bigl(\sqrt{s(\log n)/n}\bigr),\qquad\operatorname{sgn}(
\hat{\bolds{\beta}}_1)=\operatorname{sgn}\bigl(\bolds{
\beta}_1^*\bigr)\quad\mbox{and}\quad\hat{\bolds{\beta}}_2
= \mathbf{0}.
\]
If in addition, $\max_{i}\|\mathbf{H}^{1/2}\mathbf{Z}_{ni}\|_2 =
o(s^{-7/2}(\log
s)^{-1})$, then we also have
\[
\mathbf{c}^T\bigl(\mathbf{Z}_n^T
\mathbf{Z}_n\bigr)^{-1/2}\mathbf{V}_n^{-1}
\bigl(\hat{\bolds{\beta}}_1 - \bolds{\beta} _1^*\bigr)
\overset{\mathscr{D}} {\longrightarrow}N \bigl(0, \tau(1-\tau) \bigr),
\]
where $\mathbf{c}$ is an arbitrary $s$-dimensional vector satisfying
$\mathbf{c}
^T\mathbf{c}=1$.
\end{corollary}

Corollary~\ref{cor1} provides sufficient conditions for ensuring the
variable selection sign consistency of AR-Lasso. These conditions
require that R-Lasso in the initial step has the sure screening
property. We remark that in implementation, \mbox{AR-}Lasso is able to select
the variables missed by R-Lasso, as demonstrated in our numerical
studies in the next section. The theoretical comparison of the variable
selection results of R-Lasso and AR-Lasso would be an interesting topic
for future study. One set of $(p,n,s,\kappa_n)$ satisfying\vspace*{1pt} conditions
in Corollary~\ref{cor1} is $\log p= O(n^{b_1}), s = o(n^{(1-b_1)/2})$
and $\kappa_n = o(n^{b_1/4}(\log n)^{-3/2})$ with $b_1 \in(0,1/2)$
some constant.
Corollary~\ref{cor1} gives one specific choice of $\lambda_n$, not necessarily
the smallest $\lambda_n$, which makes our procedure work. In fact, the
condition on $\lambda_n$ can be weakened to $\lambda_n > 2(a+1)^{-1}\|
\hat{\bolds{\beta}}{}^{\mathrm{ini}}_1-\bolds{\beta}_1\|_\infty$.
Currently, we use the
$L_2$-norm $\|\hat{\bolds{\beta}}{}^{\mathrm{ini}}_1-\bolds{\beta}_1\|_2$
to bound this
$L_\infty$-norm, which is too crude. If one can establish $\|\hat{\bolds
{\beta}}{}^{\mathrm{ini}}_1-\bolds{\beta}_1\|_\infty= O_p(\sqrt{n^{-1}
\log p})$ for an
initial estimator $\hat{\bolds{\beta}}{}^{\mathrm{ini}}_1$, then the choice of~$\lambda_n$ can be as small as $O(\sqrt{n^{-1} \log p})$, the same
order as that used in \citet{W2012}. On the other hand, since we are
using AR-Lasso, the choice of $\lambda_n$ is not as sensitive as R-Lasso.

\section{Numerical studies}\label{secsimul}

In this section, we evaluate the finite sample property of our proposed
estimator with synthetic data. Please see the supplementary material
[\citet{ARLassoSupt}] for a real life data-set analysis, where we
provide results of an eQTL study on the \textit{CHRNA6} gene.

To assess the performance of the proposed estimator and compare it with
other methods, we simulated data from the high-dimensional linear
regression model
\[
y_{i}=\mathbf{x}_{i}^{T}\bolds{
\beta}_{0}+\varepsilon_{i},\qquad\mathbf{x}\sim\mathcal{N}
(0,\bolds{\Sigma}_{\mathbf{x}} ),
\]
where the data had $n=100$ observations and the number of parameters
was chosen as $p=400$. We fixed the true regression coefficient vector as
\[
\bolds{\beta}_0= \{ 2,0,1.5,0,0.80,0,0,1,0,1.75,0,0,0.75,0,0,0.3,0,\ldots,0 \}.
\]
For the distribution of the noise, $\varepsilon$, we considered six
symmetric distributions: normal with variance 2 ($\mathcal{N}(0,2)$),
a scale mixture of Normals for which $\sigma_i^2=1$ with probability
0.9 and $\sigma^2_i=25$ otherwise ($\mathrm{MN}_1$), a different scale mixture
model where $\varepsilon_i \sim\mathcal{N}(0,\sigma_i^2)$ and
$\sigma_i \sim\operatorname{Unif}(1,5)$ ($\mathrm{MN}_2$), Laplace, Student's $t$
with degrees of freedom 4 with doubled variance ($\sqrt{2} \times
t_4$) and Cauchy. We take $\tau=0.5$, corresponding to
$L_1$-regression, throughout the simulation. Correlation of the
covariates, $\bolds{\Sigma}_{\mathbf{x}}$ were either chosen to be identity
(i.e., $\bolds{\Sigma}_{\mathbf{x}}=I_p$) or they were generated from an
AR(1) model with correlation 0.5, that is $\bolds{\Sigma}_{\mathbf{x}_{(i,j)}}=0.5^{|i-j|}$.

We implemented five methods for each setting:
\begin{enumerate}[5.]
\item \textit{$L_2$-Oracle}, which is the least squares estimator based on
the signal covariates. 
%
\item\textit{Lasso}, the penalized least-squares estimator with
$L_1$-penalty as in \citet{tibshirani2}.
\item\textit{SCAD}, the penalized least-squares estimator with SCAD
penalty as in \citet{fan2}.
\item\textit{R-Lasso}, the robust Lasso defined as the minimizer of
(\ref{eql1reg}) with $\mathbf{d}=\bolds{1}$.
\item\textit{AR-Lasso}, which is the adaptive robust Lasso whose adaptive
weights on the penalty function were computed based on the SCAD penalty
using the R-Lasso estimate as an initial value.
\end{enumerate}

The tuning parameter, $\lambda_n$, was chosen optimally based on 100
validation data-sets. For each of these data-sets, we ran a grid search
to find the best $\lambda_n$ (with the lowest $L_2$ error for $\bolds
{\beta}
$) for the particular setting. This optimal $\lambda_n$ was recorded
for each of the 100 validation data-sets. The median of these 100
optimal $\lambda_n$ were used in the simulation studies. We preferred
this procedure over cross-validation because of the instability of the
$L_2$ loss under heavy tails.

\begin{figure}

\includegraphics{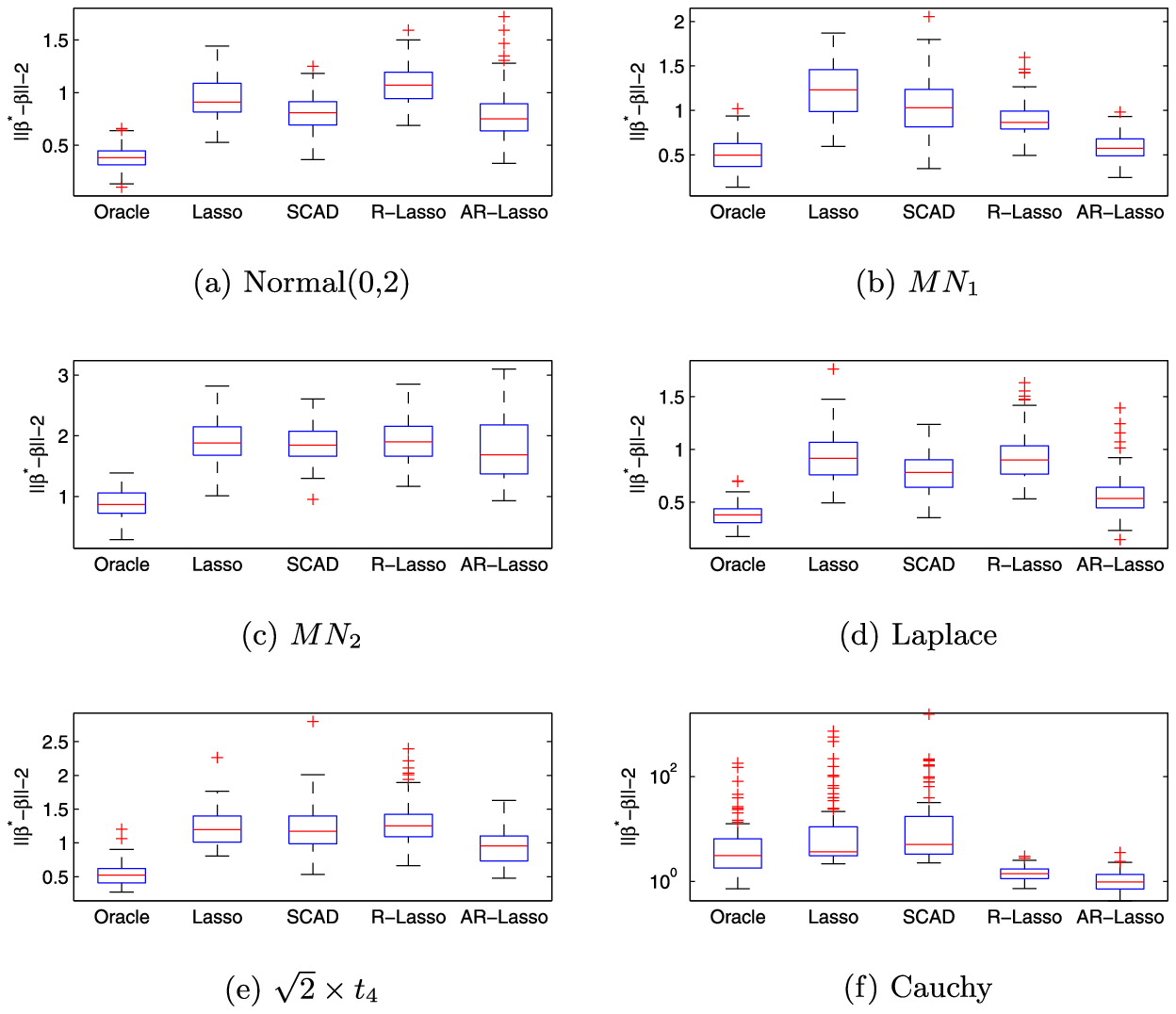}

\caption{Boxplots for $L_2$ loss with correlated covariates.}\label{figboxplot2}
\end{figure}

The following four performance measures were calculated:
\begin{enumerate}
%
\item$L_2$ loss, which is defined as $\|\bolds{\beta}^*-\hat{\bolds
{\beta}}\|_2$.
\item$L_1$ loss, which is defined as $\|\bolds{\beta}^*-\hat{\bolds
{\beta}}\|_1$.
\item Number of noise covariates that are included in the model, that
is the number of false positives (FP).
\item Number of signal covariates that are not included, that is, the
number of false negatives (FN).
\end{enumerate}
For each setting, we present the average of the performance measure
based on 100 simulations. The results are depicted in Tables~\ref
{simul1} and~\ref{simul2}. A boxplot of the $L_2$ losses under
different noise settings is also given in Figure~\ref{figboxplot2}
(the $L_2$ loss boxplot for the independent covariate setting is
similar and omitted).
For the results in Tables~\ref{simul1} and~\ref{simul2}, one should
compare the performance between Lasso and R-Lasso and that between SCAD
and AR-Lasso. This comparison reflects the effectiveness of
$L_1$-regression in dealing with heavy-tail distributions. Furthermore,
comparing Lasso with SCAD, and R-Lasso with AR-Lasso, shows the
effectiveness of using adaptive weights in the penalty function.

%
\begin{table}
\tabcolsep=0pt
\caption{Simulation results with independent covariates}\label{simul1}
\begin{tabular*}{\tablewidth}{@{\extracolsep{\fill}}@{}lcd{2.3}d{2.3}d{2.3}cc@{}}
\hline
& & \multicolumn{1}{c}{$\bolds{L_{2}}$ \textbf{Oracle}} & \multicolumn{1}{c}{\textbf{Lasso}}
  & \multicolumn{1}{c}{\textbf{SCAD}} & \multicolumn{1}{c}{\textbf{R-Lasso}} & \multicolumn{1}{c@{}}{\textbf{AR-Lasso}}
\\
\hline
$\mathcal{N}(0,2)$ & $L_{2}$ loss & 0.833 & 4.114 & 3.412 & 5.342 & 2.662\\
& $L_{1}$ loss & 0.380 & 1.047 & 0.819 & 1.169 & 0.785\\
& FP, FN & \multicolumn{1}{c}{--} & \multicolumn{1}{c}{27.00, 0.49} & \multicolumn{1}{c}{29.60, 0.51} & \multicolumn{1}{c}{36.81, 0.62} & 17.27, 0.70
\\[3pt]
$\mathrm{MN}_1$ & $L_{2}$ loss & 0.977 & 5.232 & 4.736 & 4.525 & 2.039 \\
& $L_{1}$ loss & 0.446 & 1.304 & 1.113 & 1.028 & 0.598\\
& FP, FN & \multicolumn{1}{c}{--} & \multicolumn{1}{c}{26.80, 0.73} & \multicolumn{1}{c}{29.29, 0.68} & \multicolumn{1}{c}{34.26, 0.51} & 16.76, 0.51
\\[3pt]
$\mathrm{MN}_2$ & $L_{2}$ loss & 1.886 & 7.563 & 7.583 & 8.121 & 5.647 \\
& $L_{1}$ loss & 0.861 & 2.085 & 2.007 & 2.083 & 1.845\\
& FP, FN & \multicolumn{1}{c}{--} & \multicolumn{1}{c}{20.39, 2.28} & \multicolumn{1}{c}{23.25, 2.19} & \multicolumn{1}{c}{24.64, 2.29} & 11.97, 2.57
\\[3pt]
Laplace & $L_{2}$ loss & 0.795 & 4.056 & 3.395 & 4.610 & 2.025 \\
& $L_{1}$ loss & 0.366 & 1.016 & 0.799 & 1.039 & 0.573\\
& FP, FN & \multicolumn{1}{c}{--} & \multicolumn{1}{c}{26.87, 0.62} & \multicolumn{1}{c}{29.98, 0.49} & \multicolumn{1}{c}{34.76, 0.48} & 18.81, 0.40
\\[3pt]
$\sqrt{2} \times t_4$ & $L_{2}$ loss & 1.087 & 5.303 & 5.859 & 6.185 & 3.266 \\
& $L_{1}$ loss & 0.502 & 1.378 & 1.256 & 1.403 &0.951 \\
& FP, FN & \multicolumn{1}{c}{--} & \multicolumn{1}{c}{24.61, 0.85} & \multicolumn{1}{c}{36.95, 0.76} & \multicolumn{1}{c}{33.84, 0.84} & 18.53, 0.82
\\[3pt]
Cauchy & $L_{2}$ loss & 37.451 & 211.699 & 266.088 & 6.647 & 3.587\\
& $L_{1}$ loss & 17.136 & 30.052 & 40.041 & 1.646 & 1.081\\
& FP, FN & \multicolumn{1}{c}{--} & \multicolumn{1}{c}{27.39, 5.78} & \multicolumn{1}{c}{34.32, 5.94} & \multicolumn{1}{c}{27.33, 1.41} & 17.28, 1.10\\
\hline
\end{tabular*}\vspace*{-3pt}
\end{table}

%
\begin{table}
\tabcolsep=0pt
\caption{Simulation results with correlated covariates}\label{simul2}
\begin{tabular*}{\tablewidth}{@{\extracolsep{\fill}}@{}lcd{2.3}d{3.3}d{3.3}cc@{}}
\hline
& & \multicolumn{1}{c}{$\bolds{L_{2}}$ \textbf{Oracle}} & \multicolumn{1}{c}{\textbf{Lasso}}
  & \multicolumn{1}{c}{\textbf{SCAD}} & \multicolumn{1}{c}{\textbf{R-Lasso}} & \multicolumn{1}{c@{}}{\textbf{AR-Lasso}}
\\
\hline
$\mathcal{N}(0,2)$ & $L_{2}$ loss & 0.836 & 3.440 & 3.003 & 4.185 & 2.580 \\
& $L_{1}$ loss & 0.375 & 0.943 & 0.803 & 1.079 & 0.806 \\
& FP, FN & \multicolumn{1}{c}{--} & \multicolumn{1}{c}{20.62, 0.59} & \multicolumn{1}{c}{23.13, 0.56} & \multicolumn{1}{c}{22.72, 0.77} & 14.49, 0.74
\\[3pt]
$\mathrm{MN}_1$ & $L_{2}$ loss & 1.081 & 4.415 & 3.589 & 3.652 & 1.829 \\
& $L_{1}$ loss & 0.495 & 1.211 & 1.055 & 0.901 & 0.593\\
& FP, FN & \multicolumn{1}{c}{--} & \multicolumn{1}{c}{18.66, 0.77} & \multicolumn{1}{c}{15.71, 0.75} & \multicolumn{1}{c}{26.65, 0.60} & 13.29, 0.51
\\[3pt]
$\mathrm{MN}_2$ & $L_{2}$ loss & 1.858 & 6.427 & 6.249 & 6.882 & 4.890 \\
& $L_{1}$ loss & 0.844 & 1.899 & 1.876 & 1.916 & 1.785\\
& FP, FN & \multicolumn{1}{c}{--} & \multicolumn{1}{c}{15.16, 2.08} & \multicolumn{1}{c}{14.77, 1.96} & \multicolumn{1}{c}{18.22, 1.91} & 7.86, 2.71
\\[3pt]
Laplace & $L_{2}$ loss & 0.803 & 3.341 & 2.909 & 3.606 & 1.785 \\
& $L_{1}$ loss & 0.371 & 0.931 & 0.781 & 0.927 & 0.573\\
& FP, FN & \multicolumn{1}{c}{--} & \multicolumn{1}{c}{19.32, 0.62} & \multicolumn{1}{c}{21.60, 0.38} & \multicolumn{1}{c}{24.44, 0.46} & 12.90, 0.55
\\[3pt]
$\sqrt{2} \times t_4$ & $L_{2}$ loss & 1.122 & 4.474 & 4.259 & 4.980 & 2.855 \\
& $L_{1}$ loss & 0.518 & 1.222 & 1.201 & 1.299 &0.946 \\
& FP, FN & \multicolumn{1}{c}{--} & \multicolumn{1}{c}{20.00, 0.76} & \multicolumn{1}{c}{18.49, 0.91} & \multicolumn{1}{c}{23.56, 0.79} & 13.40, 1.05
\\[3pt]
Cauchy & $L_{2}$ loss & 31.095 & 217.395 & 243.141 & 5.388 & 3.286\\
& $L_{1}$ loss & 13.978 & 31.361 & 36.624 & 1.461 & 1.074\\
& FP, FN & \multicolumn{1}{c}{--} & \multicolumn{1}{c}{25.59, 5.48} & \multicolumn{1}{c}{32.01, 5.43} & \multicolumn{1}{c}{20.80, 1.16} & 12.45, 1.17
\\
\hline
\end{tabular*}
\end{table}

%

Our simulation results reveal the following facts. The quantile based
estimators were more robust in dealing with the outliers. For example,
for the first mixture model ($\mathrm{MN}_1$) and Cauchy,
R-Lasso outperformed Lasso, and AR-Lasso outperformed SCAD in all of
the four metrics, and significantly so when the error distribution is
the Cauchy distribution. On the other hand, for the light-tail
distributions such as the normal distribution, the efficiency loss was
limited. When the tails get heavier, for instance, for the Laplace
distribution, quantile based methods started to outperform the
least-squares based approaches, more so when the tails got heavier.

The effectiveness of weights in AR-Lasso is self-evident. SCAD
outperformed Lasso and AR-Lasso outperformed R-Lasso in almost all of
the settings. Furthermore, for all of the error settings AR-Lasso had
significantly lower $L_2$ and $L_1$ loss as well as a smaller model
size compared to other estimators.

It is seen that when the noise does not have heavy tails, that is for
the normal and the Laplace distribution, all the estimators are
comparable in terms of $L_1$ loss. As expected, estimators that
minimize squared loss worked better than R-Lasso and AR-Lasso
estimators under Gaussian noise, but their performances deteriorated as
the tails got heavier. In addition, in the two heteroscedastic
settings, AR-Lasso had the best performance among others.

For Cauchy noise, least squares methods could only recover 1 or 2 of
the true variables on average. On the other hand, $L_1$-estimators
(R-Lasso and AR-Lasso) had very few false negatives, and as evident
from $L_2$ loss values, these estimators only missed variables with
smaller magnitudes. 

In addition, AR-Lasso consistently selected a smaller set of variables
than \mbox{R-}Lasso. For instance, for the setting with independent
covariates, under the Laplace distribution, R-Lasso and AR-Lasso had on
average 34.76 and 18.81 false positives, respectively. Also note that
AR-Lasso consistently outperformed R-Lasso: it estimated $\bolds{\beta}^*$
(lower $L_1$ and $L_2$ losses), and the support of $\bolds{\beta}^*$ (lower
averages for the number of false positives) more efficiently.



\section{Proofs} \label{secproof}
In this section, we prove Theorems~\ref{T1},~\ref{T2} and~\ref{T4} and provide the lemmas
used in these proofs. The proofs of Theorems~\ref{T3} and~\ref{T5} and Proposition~\ref{prop1}
are given in the supplementary Appendix [\citet{ARLassoSupt}].

We\vspace*{1pt} use techniques from empirical process theory to prove the
theoretical results. Let $v_n(\bolds{\beta}) = \sum_{i=1}^n \rho
_\tau(y_i
- \mathbf{x}_i^T\bolds{\beta})$. Then $L_n(\bolds{\beta}) =
v_n(\bolds{\beta}) +
n\lambda_n\sum_{j=1}^p d_j|\beta_j|$. For a~given deterministic
$M>0$, define the set
\[
\mathcal{B}_0(M) = \bigl\{\bolds{\beta}\in\mathbf{R}^p\dvtx
\bigl\|\bolds{\beta}- \bolds{\beta}^*\bigr\| _2\leq M, \operatorname{supp}(
\bolds{\beta})\subseteq\operatorname{supp}\bigl(\bolds{\beta}^*\bigr)
\bigr\}.
\]
Then define the function
%
%
\begin{equation}
\label{e050} Z_n(M) = \sup_{\bolds{\beta}\in\mathcal{B}_0(M)}\frac{1}{n}
\bigl|\bigl(v_n(\bolds{\beta})- v_n\bigl(\bolds{\beta}^*\bigr)
\bigr) - E \bigl(v_n(\bolds{\beta})- v_n\bigl(\bolds{
\beta} ^*\bigr) \bigr) \bigr|.
\end{equation}
Lemma~\ref{L3} in Section~\ref{lemmas} gives the rate of convergence
for $Z_n(M)$.

\subsection{\texorpdfstring{Proof of Theorem \protect\ref{T1}}{Proof of Theorem 1}}

We first show that for any $\bolds{\beta}= (\bolds{\beta
}_1^T,\mathbf{0}^T)^T\in
\mathcal{B}_0(M)$ with $M=o(\kappa_n^{-1}s^{-1/2})$,
%
%
\begin{equation}
\label{e032} E\bigl[v_n(\bolds{\beta}) - v_n\bigl(
\bolds{\beta}^*\bigr)\bigr] \geq\tfrac
{1}{2}c_0cn\bigl\|\bolds{\beta}
_1-\bolds{\beta}_1^*\bigr\|_2^2
\end{equation}
for sufficiently large $n$, where $c$ is the lower bound for $f_i(\cdot
)$ in the neighborhood of $0$. The intuition follows from the fact that
$\bolds{\beta}^*$ is the minimizer of the function $ E v_n(\bolds
{\beta})$, and
hence in Taylor's expansion of $ E[v_n(\bolds{\beta}) - v_n(\bolds
{\beta}^*)] $
around $\bolds{\beta}^*$, the first-order derivative is zero at the point
$\bolds{\beta}= \bolds{\beta}^*$. The left-hand side of (\ref
{e032}) will be
controlled by $Z_n(M)$. This yields the $L_2$-rate of convergence in
Theorem~\ref{T1}.

To prove (\ref{e032}), we set $a_i = |\mathbf{S}_i^T(\bolds{\beta
}_1-\bolds{\beta}^*_1)|$.
Then, for $\bolds{\beta}\in\mathcal{B}_0(M)$, 
\[
|a_i|\leq\|\mathbf{S}_i\|_2\bigl\|\bolds{
\beta}_1-\bolds{\beta}_1^*\bigr\| _2\leq\sqrt{s}
\kappa_nM \rightarrow0.
\]
Thus, if $\mathbf{S}_i^T(\bolds{\beta}_1-\bolds{\beta}^*_1) > 0$,
by $E1\{\varepsilon
_i\leq0\} =\tau$, Fubini's theorem, mean value theorem and Condition
\ref{assp1} it is easy to derive that
%
%
\begin{eqnarray}
\label{e014}
&& E \bigl[\rho_\tau(\varepsilon_i-a_i)-
\rho_\tau(\varepsilon_i) \bigr]\nonumber
\\
&&\qquad = E\bigl[a_i
\bigl(1\{\varepsilon_i\leq a_i\}-\tau\bigr)-
\varepsilon_i1\{ 0\leq\varepsilon_i\leq a_i
\}\bigr]
\nonumber\\[-8pt]\\[-8pt]
&&\qquad = E\biggl[\int_0^{a_i}1\{0\leq
\varepsilon_i\leq s\}\,ds\biggr]\nonumber
\\
&&\qquad = \int_0^{a_i}
\bigl(F_i(s)-F_i(0)\bigr)\,ds = \frac{1}{2}f_i(0)a_i^2
+ o(1)a_i^2,\nonumber
\end{eqnarray}
where the $o(1)$ is uniformly over all $i=1,\ldots, n$. When $\mathbf{S}
_i^T(\bolds{\beta}_1-\bolds{\beta}_1^*) < 0$, the same result can
be obtained.
Furthermore, by Condition~\ref{assp2},
\[
\sum_{i=1}^nf_i(0)a_i^2
= \bigl(\bolds{\beta}_1-\bolds{\beta}_1^*
\bigr)^T\mathbf{S}^T\mathbf{H}\mathbf{S} \bigl(\bolds{
\beta}_1-\bolds{\beta}_1^*\bigr) \geq c_0n
\bigl\|\bolds{\beta}_1-\bolds{\beta}_1^*
\bigr\|_2^2.
\]
This together with (\ref{e014}) and the definition of $v_n(\bolds
{\beta})$
proves (\ref{e032}).

The inequality (\ref{e032}) holds for any $\bolds{\beta}= (\bolds{\beta
}{}_1^T,\mathbf{0}^T)^T\in\mathcal{B}_0(M)$, yet $\hat{\bolds{\beta}}{}^o= ((\hat{\bolds{\beta}}{}_1^o)^T,\mathbf{0}^T)^T$ may not be in
the set. Thus,
we let $\tilde{\bolds{\beta}}= (\tilde{\bolds{\beta}}{}_1^T,\mathbf
{0}^T)^T$, where
\[
\tilde{\bolds{\beta}}_1 = u\hat{\bolds{\beta}}{}^o_1+(1-u)
\bolds{\beta}^*_1\qquad\mbox{with } u = M/\bigl(M + \bigl\|\hat{
\bolds{\beta}}{}_1^o -\bolds{\beta}_{1}^*
\bigr\|_2\bigr),
\]
which falls in the set $\mathcal{B}_0(M)$. Then, by the convexity and
the definition of~$\hat{\bolds{\beta}}{}_1^o$,
\[
L_n(\tilde{\bolds{\beta}}) \leq u L_n\bigl(\hat{\bolds{
\beta}}{}^o_1,\mathbf{0}\bigr) + (1-u) L_n
\bigl(\bolds{\beta}^*_1,\mathbf{0}\bigr) \leq L_n \bigl(
\bolds{\beta}^*_1,\mathbf{0} \bigr)=L_n\bigl(\bolds{
\beta}^*\bigr).
\]
Using this and the triangle inequality, we have
%
%
\begin{eqnarray}\label{eq16}
&& E\bigl[v_n(\tilde{\bolds{\beta}}) - v_n
\bigl(\bolds{\beta}^*\bigr)\bigr]\nonumber
\\
&&\qquad  =  \bigl\{ v_n\bigl(\bolds{
\beta}^*\bigr) - E v_n\bigl(\bolds{\beta}^*\bigr)\bigr\} -\bigl
\{v_n(\tilde{\bolds{\beta}}) - E v_n(\tilde{\bolds{\beta}}
)\bigr\}
\nonumber\\[-8pt]\\[-8pt]
&&\quad\qquad{} +L_n(\tilde{\bolds{\beta}}) - L_n \bigl(\bolds{
\beta}^*\bigr) + n\lambda_n\bigl\|\mathbf{d} _0\circ\bolds{
\beta}_1^*\bigr\|_1 - n\lambda_n\bigl\|
\mathbf{d}_0\circ\tilde{\bolds{\beta}} _1
\bigr\|_1\nonumber
\\
&&\qquad \leq nZ_n(M) + n\lambda_n\bigl\|\mathbf{d}_0
\circ\bigl(\bolds{\beta}_1^* -\tilde{\bolds{\beta}}_1
\bigr) \bigr\|_1.\nonumber
\end{eqnarray}
By the Cauchy--Schwarz inequality, the very last term is bounded
by\break
$n\lambda_n\|\mathbf{d}_0\|_2\|\tilde{\bolds{\beta}}_1-\bolds
{\beta}_1^*\|_2 \leq
n\lambda_n\|\mathbf{d}_0\|_2 M$.

Define the event $\mathcal{E}_n = \{Z_n(M) \leq2Mn^{-1/2}\sqrt{s\log
n}\}$. Then by Lemma~\ref{L3},
%
%
\begin{equation}
\label{e035} P(\mathcal{E}_n)\geq1- \exp\bigl({-c_0s(
\log n)/8}\bigr).
\end{equation}
On the event $\mathcal{E}_n$, by (\ref{eq16}), we have
\[
E\bigl[v_n(\tilde{\bolds{\beta}}) - v_n\bigl(\bolds{
\beta}^*\bigr)\bigr] \leq2M\sqrt{sn(\log n)}+n\lambda_n\|
\mathbf{d}_0\|_2M.
\]
Taking $M = 2\sqrt{s/n}+\lambda_n\|\mathbf{d}_0\|_2$. By Condition
\ref
{assp2} and the assumption\break  $\lambda_n\|\mathbf{d}_0\|_2\sqrt{s}\kappa
_n\rightarrow0$, it is easy to check that $M=o(\kappa
_n^{-1}s^{-1/2})$. Combining these two results with (\ref{e032}), we
obtain that on the event $\mathcal{E}_n$,
\[
\tfrac{1}{2}c_0n\bigl\|\tilde{\bolds{\beta}}_1-\bolds{
\beta}_1^*\bigr\| _2^2 \leq\bigl(2\sqrt{sn(\log
n)}+n\lambda_n\|\mathbf{d}_0\|_2\bigr) \bigl(2
\sqrt{s/n}+\lambda_n\| \mathbf{d}_0\|_2\bigr),
\]
which entails that
\[
\bigl\|\bolds{\beta}_1^*-\tilde{\bolds{\beta}}_1
\bigr\|_2 \leq O\bigl(\lambda_n\|\mathbf{d}_0
\|_2 + \sqrt{s(\log n)/n}\bigr).
\]
Note that $\|\bolds{\beta}_1^*-\tilde{\bolds{\beta}}\|_2 \leq M/2$
implies $\|
\hat{\bolds{\beta}}{}_1^o - \bolds{\beta}_1^*\|_2\leq M$. Thus,
on the event~$\mathcal{E}_n$,
\[
\bigl\|\hat{\bolds{\beta}}{}_1^o-\bolds{\beta}_1^*
\bigr\|_2\leq O\bigl(\lambda_n\| \mathbf{d}_0
\|_2 + \sqrt{s(\log n)/n}\bigr).
\]

The second result follows trivially.

\subsection{\texorpdfstring{Proof of Theorem \protect\ref{T2}}{Proof of Theorem 2}}

Since $\hat{\bolds{\beta}}{}_1^o$ defined in Theorem~\ref{T1} is a
minimizer of
$L_n(\bolds{\beta}_1, \mathbf{0})$, it satisfies the KKT conditions.
To prove that $\hat{\bolds{\beta}}=((\hat{\bolds{\beta}}{}_1^o)^T,\mathbf{0}^T)^T
\in\mathbf{R}^p$ is a global minimizer of $L_n(\bolds{\beta})$ in the
original $\mathbf{R}^p$ space, we only need to check the following condition:
%
%
\begin{equation}
\label{e005} \bigl\|\mathbf{d}_1^{-1}\circ
\mathbf{Q}^T\rho_\tau'\bigl(\mathbf{y}-
\mathbf{S}\hat{\bolds{\beta}}{}_1^o\bigr)\bigr\|
_\infty< n\lambda_n,
\end{equation}
where $\rho_\tau'(\mathbf{u})= (\rho_\tau'(u_i),\ldots, \rho
'_{\tau
}(u_n))^T$ for any $n$-vector $\mathbf{u}= (u_1,\ldots, u_n)^T$ with
$\rho
_\tau'(u_i)=\tau-1\{u_i\leq0\}$.
Here, $\mathbf{d}_1^{-1}$ denotes the vector $(d_{s+1}^{-1}, \ldots,
d_p^{-1})^T$. Then the KKT conditions and the convexity of $L_n(\bolds
{\beta}
)$ together ensure that $\hat{\bolds{\beta}}$ is a global
minimizer of
$L(\bolds{\beta})$.

Define events
\[
A_1 = \bigl\{\bigl\|\hat{\bolds{\beta}}{}_1^o -
\bolds{\beta}_1^*\bigr\|_2 \leq\gamma_n\bigr\},
\qquad A_2 = \Bigl\{\sup_{\bolds{\beta}\in\mathcal{N}}\bigl\|\mathbf
{d}_1^{-1}\circ\mathbf{Q} ^T
\rho_{\tau}'(\mathbf{y}- \mathbf{S}\bolds{
\beta}_1)\bigr\|_\infty< n\lambda_n\Bigr\},
\]
where $\gamma_n$ is defined in Theorem~\ref{T1} and
\[
\mathcal{N} = \bigl\{\bolds{\beta}= \bigl(\bolds{\beta}_1^T,
\bolds{\beta}_2^T\bigr)^T\in
R^{p}\dvtx  \bigl\| \bolds{\beta}_1 -\bolds{\beta}_1^*
\bigr\|_2 \leq\gamma_n, \bolds{\beta}_{2}=
\mathbf{0}\in\mathbf{R}^{p-s} \bigr\}.
\]
Then by Theorem~\ref{T1} and Lemma~\ref{L1} in Section~\ref{lemmas},
$ P(A_1\cap A_2) \geq1-o(n^{-cs})$.
Since $\hat{\bolds{\beta}}\in\mathcal{N}$ on the event $A_1$,
the inequality (\ref{e005}) holds on the event $A_1 \cap A_2$. This
completes the proof of Theorem~\ref{T2}.

\subsection{\texorpdfstring{Proof of Theorem \protect\ref{T4}}{Proof of Theorem 4}}

The idea of the proof follows those used in the proof of Theorems~\ref
{T1} and~\ref{T2}. We first consider the minimizer of $\widehat
L_n(\bolds{\beta})$ in the subspace $\{\bolds{\beta}= (\bolds
{\beta}_1^T,\bolds{\beta}
_2^T)^T\in\mathbf{R}^p\dvtx  \bolds{\beta}_2=\mathbf{0}\}$. Let $\bolds
{\beta}= (\bolds{\beta}
_1^T,\mathbf{0})^T$, where $\bolds{\beta}_1 = \bolds{\beta}_1^* +
\tilde a_n\mathbf{v}_1 \in\mathbf{R}^s$ with $\tilde a_n=\sqrt{s(\log
n)/n}+\lambda_n
(\|\mathbf{d}_0^*\|_2 + C_2 c_5\sqrt{s(\log p)/n} )$, $\|
\mathbf{v}_1\|_2=C$, and $C>0$ is some large enough constant. By the
assumptions in the theorem, we have $\tilde a_n = o(\kappa
_n^{-1}s^{-1/2})$. Note that
%
%
\begin{equation}
\label{e042} \widehat L_n\bigl(\bolds{\beta}^*_1 +
\tilde a_n\mathbf{v}_1,\mathbf{0}\bigr) - \widehat
L_n\bigl(\bolds{\beta}^*_1, \mathbf{0}\bigr) =
I_1(\mathbf{v}_1) + I_2(
\mathbf{v}_1),
\end{equation}
where $I_1(\mathbf{v}_1) = \|\rho_\tau(\mathbf{y}- \mathbf
{S}(\bolds{\beta}
^*_1+\tilde a_n\mathbf{v}_1))\|_1 - \|\rho_\tau(\mathbf{y}- \mathbf
{S}\bolds{\beta}
^*_1)\|_1$ and $I_2(\mathbf{v}_1) = n\lambda_n(\|\hat{\mathbf{d}
}_0\circ(\bolds{\beta}^*_1 + \tilde a_n\mathbf{v}_1)\|_1-\|\hat{\mathbf{d}
}_0\circ\bolds{\beta}^*_1\|_1)$ with\vspace*{1pt} $\|\rho_\tau(\mathbf{u})\|
_1=\sum_{i=1}^n\rho_\tau(u_i)$ for any vector $\mathbf{u}= (u_1,\ldots, u_n)^T$.
By the results in the proof of Theorem~\ref{T1}, $E[I_1(\mathbf
{v}_1)] \geq2^{-1}c_0n\|\tilde a_n\mathbf{v}_1\|_2^2$, and moreover,
with probability at least $1-n^{-cs}$,
\[
\bigl|I_1(\mathbf{v}_1)-E\bigl[I_1(
\mathbf{v}_1)\bigr]\bigr| \leq nZ_n(Ca_n)\leq2
\tilde a_n\sqrt{s(\log n)n} \|\mathbf{v}_1
\|_2.
\]
Thus, by the triangle inequality,
%
%
\begin{equation}
\label{e043} I_1(\mathbf{v}_1) \geq2^{-1}c_0
\tilde a_n^2n\|\mathbf{v}_1
\|_2^2 - 2\tilde a_n\sqrt{s(\log n)n} \|
\mathbf{v}_1\|_2.
\end{equation}
The second term on the right-hand side of (\ref{e042}) can be bounded as
%
%
\begin{equation}
\bigl|I_2(\mathbf{v}_1)\bigr| \leq n\lambda_n\bigl\|\hat{
\mathbf{d}}\circ(\tilde a_n\mathbf{v}_1)
\bigr\|_1 \leq na_n\lambda_n\|\hat{
\mathbf{d}}_0\| _2\| \mathbf{v}_1
\|_2.
\end{equation}
By triangle inequality and Conditions~\ref{assp5} and~\ref{assp4}, it
holds that
%
%
\begin{eqnarray}\label{e044}
\|\hat{\mathbf{d}}_0\|_2 &\leq& \bigl\|\hat{
\mathbf{d}}_0-\mathbf{d}^*_0\bigr\|_2 + \bigl\|
\mathbf{d} _0^*\bigr\|_2
\nonumber\\[-8pt]\\[-8pt]
& \leq& c_5\bigl\|\hat{\bolds{
\beta}}{}^{\mathrm{ini}}_1-\bolds{\beta}^*_1
\bigr\|_2 + \bigl\| \mathbf{d}_0^*\bigr\|_2 \leq
C_2 c_5\sqrt{s(\log p)/n}+\bigl\|\mathbf{d}_0^*
\bigr\|_2.\nonumber
\end{eqnarray}
Thus, combining (\ref{e042})--(\ref{e044}) yields
\begin{eqnarray*}
\widehat L_n\bigl(\bolds{\beta}^* + \tilde a_n
\mathbf{v}_1\bigr) - \widehat L_n\bigl(\bolds{\beta}^*
\bigr) &\geq& 2^{-1}c_0na_n^2\|
\mathbf{v}_1\|_2^2 - 2\tilde a_n
\sqrt{s(\log n)n} \|\mathbf{v}_1\|_2
\\
&&{}-na_n\lambda_n \bigl(\bigl\|\mathbf{d}_0^*
\bigr\|_2+C_2 c_5\sqrt{s(\log p)/n} \bigr)\|
\mathbf{v}_1\|_2.
\end{eqnarray*}
Making $\|\mathbf{v}_1\|_2=C$ large enough, we obtain that with
probability tending to one, $\widehat L_n(\bolds{\beta}^* + \tilde
a_n\mathbf{v}) - \widehat L_n(\bolds{\beta}^*)>0$. Then it follows
immediately that
with asymptotic probability one, there exists a minimizer $\hat{\bolds
{\beta}}_1$ of $\widehat L_n(\bolds{\beta}_1, \mathbf{0})$
such that $\|\hat{\bolds{\beta}}_1 - \bolds{\beta}^*_1\|_2 \leq C_3
\tilde a_n \equiv
a_n$ with some
constant $C_3>0$.

It remains to prove that with asymptotic probability one,
%
%
\begin{equation}
\label{e045} \bigl\|\hat{\mathbf{d}}_1^{-1}\circ
\mathbf{Q}^T\rho_\tau'(\mathbf{y}- \mathbf{S}
\hat{\bolds{\beta}}_1)\bigr\|_\infty< n\lambda_n.
\end{equation}
Then by KKT conditions, $\hat{\bolds{\beta}}= (\hat{\bolds{\beta}}{}_1^T,
\mathbf{0}^T)^T$ is a global minimizer of $\widehat L_n(\bolds{\beta})$.

Now we proceed to prove (\ref{e045}). Since $\beta^*_j=0$ for all
$j=s+1, \ldots, p$, we have that $d^*_j = p'_{\lambda_n}(0+)$.
Furthermore, by Condition~\ref{assp5}, it holds that $|\hat\beta
_j^{\mathrm{ini}}|\leq C_2\sqrt{s(\log p)/n}$ with asymptotic probability one.
Then, it follows that
\[
\min_{j>s}p'_{\lambda_n}\bigl(\bigl|\hat
\beta_j^{\mathrm{ini}}\bigr|\bigr) \geq p'_{\lambda
_n}
\bigl(C_2\sqrt{s(\log p)/n}\bigr).
\]
Therefore, by Condition~\ref{assp4} we conclude that
%
%
\begin{equation}
\label{e048} \bigl\|(\hat{\mathbf{d}}_1)^{-1}
\bigr\|_{\infty} = \Bigl(\min_{j>s}p'_{\lambda
_n}
\bigl(\bigl|\hat\beta_j^{\mathrm{ini}}\bigr|\bigr)\Bigr)^{-1}<
2/p'_{\lambda_n}(0+) = 2\bigl\| \bigl(\mathbf{d} ^*_1
\bigr)^{-1}\bigr\|_{\infty}.
\end{equation}

From the conditions of Theorem~\ref{T2} with $\gamma_n=a_n$, it
follows from Lemma~\ref{L1} [inequality (\ref{e049})] that, with
probability at least $1-o(p^{-c})$,
%
%
\begin{equation}
\label{e046} \sup_{\|\bolds{\beta}_1 -\bolds{\beta}_1^*\|_2\leq C_3a_n}
\bigl\| \mathbf{Q}^T
\rho_\tau'(\mathbf{y}- \mathbf{S}\bolds{
\beta}_1)\bigr\|_\infty< \frac{n
\lambda_n}{2\|(\mathbf{d}
_1^*)^{-1}\|_\infty}\bigl(1+o(1)\bigr).
\end{equation}
Combining (\ref{e048})--(\ref{e046}) and by the triangle inequality,
it holds that with asymptotic probability one,
\[
\sup_{\|\bolds{\beta}_1 -\bolds{\beta}_1^*\|_2\leq C_3a_n} \bigl\| (\hat
{\mathbf{d} }_1)^{-1}
\circ\mathbf{Q}^T\rho_\tau'(\mathbf{y}-
\mathbf{S}\bolds{\beta}_1)\bigr\|_\infty< n
\lambda_n.
\]
Since the minimizer $\hat{\bolds{\beta}}_1$ satisfies $\|\hat{\bolds
{\beta}}
_1 - \bolds{\beta}_1^*\|_2< C_3a_n$ with asymptotic probability one, the
above inequality ensures that (\ref{e045}) holds with probability
tending to one. This completes the proof.

\subsection{Lemmas}\label{lemmas}

This subsection contains lemmas used in proofs of Theorems~\ref{T1},~\ref{T2} and~\ref{T4}.
%
%
\begin{lemma}\label{L3} Under Condition~\ref{assp2}, for any $t>0$,
we have
%
%
\begin{equation}
P\bigl(Z_n(M) \geq4M\sqrt{s/n} + t\bigr) \leq\exp
\bigl(-nc_0t^2/\bigl(8M^2\bigr) \bigr).
\end{equation}
\end{lemma}
\begin{pf} Define $\rho(s,y)=(y-s)(\tau- 1\{y-s\leq0\})$. Then
$v_n(\bolds{\beta})$ in (\ref{e050}) can be rewritten as
$v_n(\bolds{\beta}) = \sum_{i=1}^n\rho(\mathbf{x}_i^T\bolds{\beta
},y_i)$. Note
that the following Lipschitz condition holds for $\rho(\cdot, y_i)$:
%
%
\begin{equation}
\label{e001} \bigl|\rho(s_1, y_i) - \rho(s_2,y_i)
\bigr|\leq\max\{\tau,1-\tau\} |s_1-s_2| \leq|s_1-s_2|.
\end{equation}
Let $W_1, \ldots, W_n$ be a Rademacher sequence, independent of model
errors $\varepsilon_1, \ldots, \varepsilon_n$. The Lipschitz
inequality (\ref{e001}) combined with the symmetrization
theorem and concentration inequality [see, e.g., Theorems 14.3 and 14.4
in \citet{BV11}] yields that
%
%
\begin{eqnarray}\label{e033}
E\bigl[Z_n(M)\bigr]&\leq& 2 E\sup
_{\bolds{\beta}\in\mathcal
{B}_0(M)} \Biggl|\frac{1}{n}\sum_{i=1}^n
W_i \bigl(\rho\bigl(\mathbf{x}_i^T\bolds{
\beta}, y_i\bigr) - \rho\bigl(\mathbf{x}_i^T
\bolds{\beta}^*, y_i\bigr) \bigr) \Biggr|
\nonumber\\[-8pt]\\[-8pt]
&\leq& 4E\sup_{\bolds{\beta}\in\mathcal{B}_0(M)} \Biggl|\frac
{1}{n}\sum
_{i=1}^n W_i \bigl(
\mathbf{x}_i^T\bolds{\beta}- \mathbf{x}_i^T
\bolds{\beta}^* \bigr) \Biggr|.\nonumber
\end{eqnarray}
On the other hand, by the Cauchy--Schwarz inequality
\begin{eqnarray*}
\Biggl\llvert\sum_{i=1}^n
W_i \bigl(\mathbf{x}_i^T\bolds{\beta}-
\mathbf{x}_i^T\bolds{\beta}^* \bigr)\Biggr\rrvert&=&
\Biggl\llvert\sum_{j=1}^s \Biggl(\sum
_{i=1}^nW_ix_{ij}
\Biggr) \bigl(\beta_j-\beta_j^*\bigr)\Biggr\rrvert
\\
&\leq& \bigl\|\bolds{\beta} _1-\bolds{\beta}^*_1
\bigr\|_2 \Biggl\{\sum_{j=1}^s \Biggl|\sum
_{i=1}^nW_ix_{ij}
\Biggr|^2 \Biggr\}^{1/2}.
\end{eqnarray*}
By Jensen's inequality and concavity of the square root function, $E
(X^{1/2} )\leq(E X )^{1/2}$ for any
nonnegative random variable $X$. Thus, these two inequalities ensure
that the very right-hand side of (\ref{e033}) can be further bounded by
%
%
\begin{eqnarray}\label{e034}
&& \sup_{\bolds{\beta}\in\mathcal{B}_0(M)}\bigl\|\bolds{\beta}- \bolds{
\beta}^*\bigr\| _2E \Biggl\{\sum_{j=1}^s
\Biggl|\frac{1}{n}\sum_{i=1}^nW_ix_{ij}
\Biggr|^2 \Biggr\}^{1/2}
\nonumber\\[-8pt]\\[-8pt]
&&\qquad \leq M \Biggl\{\sum_{j=1}^s E \Biggl|
\frac{1}{n}\sum_{i=1}^nW_ix_{ij}
\Biggr|^2 \Biggr\}^{1/2}= M\sqrt{s/n}.\nonumber
\end{eqnarray}
Therefore, it follows from (\ref{e033}) and (\ref{e034}) that
%
%
\begin{equation}
\label{e051} E \bigl[Z_n(M)\bigr] \leq4M\sqrt{s/n}.
\end{equation}

Next, since $n^{-1}\mathbf{S}^T\mathbf{S}$ has bounded eigenvalues,
for any $\bolds{\beta}
=(\bolds{\beta}_1^T, \mathbf{0}^T) \in\mathcal{B}_0(M)$,
\[
\frac{1}{n}\sum_{i=1}^n\bigl(
\mathbf{x}_i^T\bigl(\bolds{\beta}- \bolds{\beta}^*\bigr)
\bigr)^2 = \frac{1}{n}\bigl(\bolds{\beta}_1 -
\bolds{\beta}^*_1\bigr)^T\mathbf{S}^T
\mathbf{S}\bigl(\bolds{\beta}_1 - \bolds{\beta} ^*_1
\bigr)\leq c_0^{-1}\|\bolds{\beta}_1 -\bolds{
\beta}_1^*\|_2^2 \leq c_0^{-1}M^2.
\]
Combining this with the Lipschitz inequality (\ref{e001}), (\ref
{e051}) and applying
Massart's concentration theorem [see Theorem 14.2 in \citet{BV11}]
yields that for any $t > 0$,
\[
P\bigl(Z_n(M) \geq4M\sqrt{s/n} + t\bigr) \leq\exp
\bigl(-nc_0t^2/\bigl(8M^2\bigr) \bigr).
\]
This proves the lemma.
\end{pf}

%
\begin{lemma}\label{L1} Consider a ball in $R^{s}$ around $\bolds
{\beta}^*\dvtx \mathcal{N} = \{\bolds{\beta}= (\bolds{\beta}_1^T, \bolds{\beta
}_2^T)^T\in R^{p}\dvtx
\bolds{\beta}_{2}=0, \|\bolds{\beta}_1 -\bolds{\beta}_1^*\|_2
\leq\gamma_n\}$ with
some sequence $\gamma_n \rightarrow0$. Assume that\break  $\min_{j>s}\,d_j
>c_3$, $\sqrt{1+\gamma_ns^{3/2}\kappa_n^2}\log_2n = o(\sqrt{n} \lambda
_n)$, $n^{1/2}\lambda_n(\log p)^{-1/2}\rightarrow\infty$,
and $\kappa_n\gamma_n^2 = o(\lambda_n)$. 
Then under Conditions~\ref{assp1}--\ref{assp3}, there exists some
constant $c>0$ such that
\[
P \Bigl(\sup_{\bolds{\beta}\in\mathcal{N}}\bigl\|\mathbf{d}_1^{-1}
\circ\mathbf{Q}^T\rho_\tau'(\mathbf{y}-
\mathbf{S}\bolds{\beta}_1)\bigr\|_\infty\geq n
\lambda_n \Bigr) \leq o\bigl(p^{-c}\bigr),
\]
where $\rho_\tau'(u) = \tau-1\{u\leq0\}$.
\end{lemma}

\begin{pf}
For a fixed $j \in\{s+1, \ldots, p\}$ and $\bolds{\beta}= (\bolds
{\beta}_1^T,
\bolds{\beta}_2^T)^T\in\mathcal{N}$, define
\[
\gamma_{\bolds{\beta}, j}(\mathbf{x}_i,y_i) =
x_{ij} \bigl[\rho_\tau'\bigl(y_i-
\mathbf{x}_i^T\bolds{\beta}\bigr)-\rho_\tau'(
\varepsilon_i)-E\bigl[\rho_\tau'
\bigl(y_i-\mathbf{x}_i^T\bolds{\beta}\bigr)-
\rho_\tau'(\varepsilon_i)\bigr] \bigr],
\]
where $\mathbf{x}_i^T= (x_{i1}, \ldots, x_{ip})$ is the $i$th row of
the design matrix. The key for the proof is to use the following decomposition:
%
%
\begin{eqnarray}\label{e039}
&& \sup_{\bolds{\beta}\in\mathcal{N}} \biggl\|\frac
{1}{n}
\mathbf{Q}^T\rho_\tau'(\mathbf{y}- \mathbf{S}
\bolds{\beta}_1) \biggr\|_\infty\nonumber
\\
&&\qquad \leq  \sup_{\bolds{\beta}\in
\mathcal{N}}
\biggl\|\frac{1}{n}\mathbf{Q}^TE\bigl[\rho_\tau'(
\mathbf{y}- \mathbf{S}\bolds{\beta} _1)-\rho_\tau'(
\bolds{\varepsilon})\bigr] \biggr\|_\infty
\\
&&\quad\qquad {}+ \biggl\|\frac{1}{n}\mathbf{Q}^T\rho_\tau'(
\bolds{\varepsilon}) \biggr\|_\infty+\max_{j>s}\sup
_{\bolds{\beta}\in\mathcal{N}}\frac{1}{n}\sum_{i=1}^n
\bigl|\gamma_{\bolds{\beta},j}(\mathbf{x}_i,y_i) \bigr|.\nonumber
\end{eqnarray}
We will prove that with probability at least $1-o(p^{-c})$,
%
%
\begin{eqnarray}
\qquad I_1 &\equiv&\sup_{\bolds{\beta}\in\mathcal{N}} \biggl\|\frac
{1}{n}
\mathbf{Q} ^TE\bigl[\rho_\tau'(\mathbf{y}-
\mathbf{S}\bolds{\beta}_1)-\rho_\tau'(
\bolds{\varepsilon})\bigr] \biggr\| _\infty< \frac{\lambda_n}{2\|\mathbf
{d}_1^{-1}\|_\infty}+ o(\lambda
_n),\label{e049}
\\
I_2&\equiv& n^{-1}\bigl\|\mathbf{Q}^T
\rho_\tau'(\bolds{\varepsilon})\bigr\| _\infty= o(
\lambda_n),\label{e013}
\\
I_3 &\equiv&\max_{j>s}\sup_{\bolds{\beta}\in\mathcal{N}}
\Biggl|\frac
{1}{n}\sum_{i=1}^n
\gamma_{\bolds{\beta},j}(\mathbf{x}_i,y_i) \Biggr| =
o_p(\lambda_n).\label{e038}
\end{eqnarray}
Combining (\ref{e039})--(\ref{e038}) with the assumption $\min_{j>s}\,d_j
>c_3$ completes the proof of the lemma.

Now we proceed to prove (\ref{e049}). Note that $I_1$ can be rewritten as
%
%
\begin{equation}
\label{e019} I_1=\max_{j>s}\sup
_{\bolds{\beta}\in\mathcal{N}} \Biggl|\frac
{1}{n}\sum_{i=1}^nx_{ij}E
\bigl[\rho_\tau'(\varepsilon_i)-
\rho_\tau'\bigl(y_i-\mathbf
{x}_i^T\bolds{\beta}\bigr)\bigr] \Biggr|.
\end{equation}
By Condition~\ref{assp1},
\begin{eqnarray*}
&& E \bigl[\rho_\tau'(\varepsilon_i)-
\rho_\tau'\bigl(y_i-\mathbf{x}_i^T
\bolds{\beta} \bigr)\bigr]
\\
&&\qquad = F_i\bigl(\mathbf{S}_i^T
\bigl(\bolds{\beta}_1-\bolds{\beta}_1^*\bigr)\bigr) -
F_i(0)= f_i(0)\mathbf{S} _i^T
\bigl(\bolds{\beta}_1 - \bolds{\beta}_1^*\bigr)+
\widetilde I_{i},
\end{eqnarray*}
where $F(t)$ is the cumulative distribution function of $\varepsilon
_i$, and
$
\widetilde I_{i} = F_i(\mathbf{S}_i^T(\bolds{\beta}_{1}-\bolds{\beta
}_1^*))-F_i(0) -
f_i(0)\mathbf{S}_i^T(\bolds{\beta}_{1}-\bolds{\beta}_1^*)
$.
Thus, for any $j>s$,
\begin{eqnarray*}
 \sum_{i=1}^nx_{ij} E \bigl[
\rho_\tau'(\varepsilon_i)-
\rho_\tau'\bigl(y_i-\mathbf{x}_i^T
\bolds{\beta}\bigr)\bigr]&=& \sum_{i=1}^n
\bigl(f_i(0)x_{ij}\mathbf{S} _i^T
\bigr) \bigl(\bolds{\beta}_1-\bolds{\beta}_1^*\bigr) +
\sum_{i=1}^n x_{ij}\widetilde
I_{i}.
\end{eqnarray*}
This together with (\ref{e019}) and Cauchy--Schwarz inequality
entails that
%
%
\begin{equation}
\label{e016} I_1\leq\biggl\|\frac{1}{n}\mathbf{Q}^T
\mathbf{H}\mathbf{S}\bigl(\bolds{\beta}_1-\bolds{
\beta}_1^*\bigr)\biggr\|_{\infty} + \max_{j>s}\Biggl|
\frac{1}{n}\sum_{i=1}^n
x_{ij}\widetilde I_{i}\Biggr|,
\end{equation}
where $\mathbf{H}= \operatorname{diag}\{f_1(0,\ldots, f_n(0))\}$. We
consider the
two terms on the right-hand side of (\ref{e016}) one by one. By
Condition~\ref{assp3}, the first term can be bounded as
%
%
\begin{equation}
\label{e017} \qquad \biggl\|\frac{1}{n}\mathbf{Q}^T\mathbf{H}\mathbf{S}
\bigl(\bolds{\beta}_1-\bolds{\beta}_1^*\bigr)
\biggr\|_{\infty} \leq\biggl\| \frac{1}{n}\mathbf{Q}^T\mathbf{H}
\mathbf{S}\biggr\|_{2,\infty}\bigl\|\bolds{\beta}_1-\bolds{
\beta}_1^*\bigr\|_2 < \frac{\lambda_n}{2\|\mathbf{d}_1^{-1}\|_\infty}.
\end{equation}
By Condition~\ref{assp1}, $|\widetilde I_{i}| \leq c(\mathbf
{S}_i^T(\bolds{\beta}
_{1}-\bolds{\beta}_1^*))^2$. This together with Condition~\ref{assp2}
ensures that the second term of (\ref{e016}) can be bounded as
\begin{eqnarray*}
\max_{j>s} \Biggl| \frac{1}{n}\sum
_{i=1}^n x_{ij}\widetilde I_{i}
\Biggr| &\leq& \frac{\kappa_n}{n}\sum_{i=1}^n |
\widetilde I_{1} |
\\
&\leq& C\frac{\kappa_n}{n}\sum
_{i=1}^n\bigl(\mathbf{S}_i^T
\bigl(\bolds{\beta}_1 - \bolds{\beta}_1^*\bigr)
\bigr)^2
\\
&\leq& C\kappa_n\bigl\|\bolds{\beta}_1-
\bolds{\beta}_1^*\bigr\|_2^2.
\end{eqnarray*}
Since $\bolds{\beta}\in\mathcal{N}$, it follows from the assumption
$\lambda_n^{-1}\kappa_n\gamma_n^2 = o(1)$ that
\[
\max_{j>s} \Biggl|\frac{1}{n}\sum
_{i=1}^n x_{ij}\widetilde I_{i}
\Biggr| \leq C\kappa_n\gamma_n^2 = o(
\lambda_n).
\]
Plugging the above inequality and (\ref{e017}) into (\ref{e016})
completes the proof of (\ref{e049}).

Next, we prove (\ref{e013}). By Hoeffding's inequality, if $\lambda
_n>2\sqrt{(1+c)(\log p)/n}$ with $c$ is some positive constant, then
\begin{eqnarray*}
&& P \bigl(\bigl\|\mathbf{Q}^T\rho_\tau'(\bolds{
\varepsilon})\bigr\|_\infty\geq n\lambda_n\bigr)
\\
&&\qquad \leq \sum
_{j=s+1}^{p}2\exp\biggl(-\frac{n^2\lambda_n^2}{4\sum_{i=1}^nx_{ij}^2}
\biggr)
 = 2\exp\bigl(\log(p-s)- n\lambda_n^2/4 \bigr)\leq O
\bigl(p^{-c}\bigr).
\end{eqnarray*}
Thus, with probability at least $1-O(p^{-c})$, (\ref{e013}) holds.

We now apply Corollary 14.4 in \citet{BV11} to prove~(\ref{e038}). To
this end, we need to check conditions of the corollary. For each fixed
$j$, define the functional space $\Gamma_j = \{\gamma_{\bolds{\beta},j}\dvtx
\bolds{\beta}\in\mathcal{N}\}$. First note that $E[\gamma_{\bolds
{\beta},
j}(\mathbf{x}_i,y_i)]=0$ for any $\gamma_{\bolds{\beta},j}\in
\Gamma_j$.
Second, since the $\rho'_\tau$ function is bounded, we have
\begin{eqnarray*}
&& \frac{1}{n}\sum_{i=1}^n
\gamma_{\bolds{\beta}, j}^2(\mathbf{x}_i, y_i)
\\
&&\qquad =
\frac{1}{n}\sum_{i=1}^nx_{ij}^2
\bigl(\rho'_\tau\bigl(y_i-\mathbf
{x}_i^T\bolds{\beta}\bigr)-\rho'_\tau(
\varepsilon_i)
 - E \bigl(\rho'_\tau\bigl(y_i-
\mathbf{x}_i^T\bolds{\beta}\bigr)-\rho_\tau(
\varepsilon_i) \bigr) \bigr)^2 \leq4.
\end{eqnarray*}
Thus, $\|\gamma_{\bolds{\beta},j}\|_n \equiv(n^{-1}\sum_{i=1}^n\gamma
_{\bolds{\beta}, j}^2(\mathbf{x}_i, y_i)^2
)^{1/2}\leq2$.

Third, we will calculate the covering number of the functional space
$\Gamma_j$, $N(\cdot,\Gamma_j,\|\cdot\|_2)$. For any $\bolds{\beta}=
(\bolds{\beta}_1^T,\bolds{\beta}_2^T)^T\in\mathcal{N}$ and
$\tilde{\bolds{\beta}}=
(\tilde{\bolds{\beta}}{}_1^T,\tilde{\bolds{\beta}}{}_2^T)^T\in
\mathcal{N}$, by
Condition~\ref{assp1} and the mean value theorem,
%
%
\begin{eqnarray}\label{e036}
&& E \bigl[\rho_\tau'\bigl(y_i-
\mathbf{x}_i^T\bolds{\beta}\bigr)-\rho_\tau
'(\varepsilon_i)\bigr]-E \bigl[\rho_\tau'
\bigl(y_i-\mathbf{x}_i^T\tilde{\bolds{
\beta}}\bigr)-\rho_\tau'(\varepsilon_i)
\bigr]\nonumber
\\
&&\qquad = F_i\bigl(\mathbf{S}_i^T\bigl(
\tilde{\bolds{\beta}}_{1}-\bolds{\beta}_1^*\bigr)
\bigr)-F_i\bigl(\mathbf{S} _i^T\bigl(\bolds{
\beta}_{1}-\bolds{\beta}_1^*\bigr)\bigr)
\\
&&\qquad = f_i(a_{1i})\mathbf{S}_i^T(
\bolds{\beta}_1 - \tilde{\bolds{\beta}}_1),\nonumber
\end{eqnarray}
where $F(t)$ is the cumulative distribution function of $\varepsilon
_i$, and $a_{1i}$ lies on the segment connecting $\mathbf
{S}_i^T(\tilde{\bolds{\beta}}_{1}-\bolds{\beta}_1^*)$ and $\mathbf
{S}_i^T(\bolds
{\beta}_{1}-\bolds{\beta}_1^*)$. Let
$\kappa_n = \max_{ij}|x_{ij}|$. Since $f_i(u)$'s are uniformly
bounded, by (\ref{e036}),
%
%
\begin{eqnarray}\label{e018}
&& \bigl|x_{ij}E \bigl[\rho_\tau'
\bigl(y_i-\mathbf{x}_i^T\bolds{\beta}\bigr)-
\rho_\tau'(\varepsilon_i)
\bigr]-x_{ij}E \bigl[\rho_\tau'
\bigl(y_i-\mathbf{x}_i^T\tilde{\bolds{
\beta}}\bigr)-\rho_\tau'(\varepsilon_i)\bigr]\bigr|\nonumber
\\
&&\qquad \leq C\bigl|x_{ij}\mathbf{S}_i^T(\bolds{
\beta}_1-\tilde{\bolds{\beta}}_1)\bigr|
\\
&&\qquad \leq C
\|x_{ij}\mathbf{S}_i\|_2\|\bolds{
\beta}_1-\tilde{\bolds{\beta}}_1\|_2 \leq C
\sqrt{s}\kappa_n^2\|\bolds{\beta}_1-\tilde{
\bolds{\beta}}_1\|_2,\nonumber
\end{eqnarray}
where $C>0$ is some generic constant. It is known [see, e.g., Lemma
14.27 in \citet{BV11}] that the ball $\mathcal{N}$ in $\mathbf{R}^s$
can be covered by $(1+4\gamma_n/\delta)^s$ balls with radius $\delta
$. 
Since $\rho_\tau'(y_i-\mathbf{x}_i^T\bolds{\beta})-\rho_\tau
'(\varepsilon_i)$ can only take 3 different values $\{-1, 0,1\}$, it
follows from (\ref{e018}) that\vspace*{1pt} the covering number of $\Gamma_j$ is
$N(2^{2-k},\Gamma_j,\|\cdot\|_2) = 3(1+C^{-1}2^{k}\gamma
_ns^{1/2}\kappa_n^2)^s$. Thus, by calculus, for any $0\leq k\leq(\log
_2 n)/2$,
\begin{eqnarray*}
&& \log\bigl(1+N\bigl(2^{2-k},\Gamma,\|\cdot\|_2\bigr)
\bigr)
\\
&&\qquad \leq\log(6) +s\log\bigl(1+C^{-1}2^{k}
\gamma_ns^{1/2}\kappa_n^2 \bigr)
\\
&&\qquad \leq\log(6) + C^{-1}2^{k}\gamma_ns^{3/2}
\kappa_n^2\leq4 \bigl(1+C^{-1}
\gamma_ns^{3/2}\kappa_n^2
\bigr)2^{2k}.
\end{eqnarray*}
Hence, conditions of Corollary 14.4 in \citet{BV11} are checked and we
obtain that for any $t>0$,
\begin{eqnarray*}
&& P \Biggl(\sup_{\bolds{\beta}\in\mathcal{N}} \Biggl|\frac{1}{n}\sum
_{i=1}^n\gamma_{\bolds{\beta},j}(\mathbf{x}_i,y_i)
\Biggr|\geq\frac
{8}{\sqrt{n}} \Bigl(3\sqrt{1+C^{-1}
\gamma_ns^{3/2}\kappa_n^2}\log
_2n + 4 +4t \Bigr) \Biggr)
\\
&&\qquad \leq4\exp\biggl(-\frac{nt^2}{8}\biggr).
\end{eqnarray*}
Taking $t = \sqrt{C(\log p)/n}$ with $C>0$ large enough constant, we
obtain that
\begin{eqnarray*}
&& P \Biggl(\max_{j>s}\sup_{\bolds{\beta}\in\mathcal{N}} \Biggl|
\frac
{1}{n}\sum_{i=1}^n
\gamma_{\bolds{\beta},j}(\mathbf{x}_i,y_i) \Biggr|  \geq
\frac{24}{\sqrt{n}}\sqrt{1+C^{-1}\gamma_ns^{3/2}
\kappa_n^2}\log_2n \Biggr)
\\
&&\qquad \leq 4(p-s)\exp\biggl(-\frac{C\log p}{8} \biggr)\rightarrow0.
\end{eqnarray*}
Thus, if $\sqrt{1+\gamma_ns^{3/2}\kappa_n^2}\log_2n = o(\sqrt{n} \lambda
_n)$, then with probability at least $1-o(p^{-c})$, (\ref{e038}) holds. This completes the proof of the lemma.
\end{pf}

\section*{Acknowledgments}
The authors sincerely thank the Editor, Associate Editor, and three
referees for their constructive comments that led to substantial
improvement of the paper.

\begin{supplement}
\stitle{Supplementary material for: Adaptive robust variable selection}
\slink[doi]{10.1214/13-AOS1191SUPP} 
\sdatatype{.pdf}
\sfilename{AOS1191\_supp.pdf}
\sdescription{Due to space constraints, the proofs of Theorems \ref{T3} and \ref{T5}
and the results of the real life data-set study are relegated to the
supplement [\citet{ARLassoSupt}].}
\end{supplement}



%

\printaddresses

\end{document}